\documentclass{amsart}
\usepackage{amssymb}
\usepackage[english]{babel}
\usepackage{amsmath}
\usepackage{latexsym}
\usepackage{eepic}
\usepackage{epsfig}
\usepackage{graphicx}
\usepackage{pb-diagram,lamsarrow,pb-lams,amscd}
\usepackage{eufrak}
\usepackage{mathrsfs}

\newtheorem{e-proposition}[theorem]{Proposition}

\newtheorem{e-definition}[theorem]{Definition\rm}

\newtheorem{example}{\it Example\/}

\def\A{\mathbb{A}}

\newcommand{\N}{\ensuremath{\mathbb{N}}}
\newcommand{\Z}{\ensuremath{\mathbb{Z}}}
\newcommand{\Q}{\ensuremath{\mathbb{Q}}}

\newcommand{\R}{\ensuremath{\mathbb{R}}}
\newcommand{\C}{\ensuremath{\mathbb{C}}}

\newcommand{\cA}{\ensuremath{\mathcal{A}}}

\newcommand{\cM}{\ensuremath{\mathcal{M}}}
\newcommand{\cH}{\ensuremath{\mathscr{H}}}

\renewcommand{\R}{\ensuremath{\mathbb{R}}}
\renewcommand{\C}{\ensuremath{\mathbb{C}}}

\newcommand{\mU}{\ensuremath{\mathfrak{U}}}

\renewcommand{\cA}{\ensuremath{\mathscr{A}}}
\renewcommand{\cM}{\ensuremath{\mathscr{M}}}
\renewcommand{\cH}{\ensuremath{\mathscr{H}}}

\newcommand{\Spec}{\ensuremath{\mathrm{Spec}\,}}
\newcommand{\Spf}{\ensuremath{\mathrm{Spf}\,}}
\newcommand{\Sp}{\ensuremath{\mathrm{Sp}\,}}

\newcommand{\Proj}{\ensuremath{\mathbb{P}}}

\begin{document}
\title[Formal and rigid geometry, and some applications]{Formal and rigid geometry: an intuitive introduction, and some applications}
\author{Johannes Nicaise}
\address{Universit\'e Lille 1\\
Laboratoire Painlev\'e, CNRS - UMR 8524\\ Cit\'e Scientifique\\59 655 Villeneuve d'Ascq C\'edex \\
France} \email{johannes.nicaise@math.univ-lille1.fr}
\begin{abstract}
We give an informal introduction to formal and rigid geometry, and
we discuss some applications in algebraic and arithmetic geometry
and singularity theory, with special emphasis on recent
applications to the Milnor fibration and the motivic zeta function
by J. Sebag and the author.
\end{abstract}
 \maketitle

\section{Introduction}
Let $R$ be a complete discrete valuation ring, with quotient field
$K$, and residue field $k$. We choose a uniformizing parameter
$\pi$, i.e. $\pi$ generates the unique maximal ideal of $R$.
Geometers may take $R=\C[[t]]$, the ring of formal power series
over the complex numbers, with $K=\C((t)),\,k=\C,\,\pi=t$, while
number theorists might prefer to think of $R=\Z_p$, the ring of
$p$-adic integers, with $K=\Q_p,\,k=\mathbb{F}_p,\,\pi=p$.

Very roughly, a formal scheme over $R$ consists of an algebraic
variety over $k$,
 together with algebraic information on an
infinitesimal neighbourhood of this variety. If $X$ is a variety
over $R$, we can associate to $X$ its formal completion
$\widehat{X}$ in a natural way. It is a formal scheme over $R$,
and can be seen as an infinitesimal tubular neighbourhood of the
special fiber $X_0$ in $X$. Its underlying topological space
coincides with the space underlying $X_0$, but additional
infinitesimal information is contained in the sheaf of regular
functions on $\widehat{X}$.

An important aspect of the formal scheme $\widehat{X}$ is the
following phenomenon. A closed point $x$ on the scheme-theoretic
generic fiber of $X$ over $K$ has coordinates in some finite
extension of the field $K$, and (unless $X$ is proper over $R$)
there is no natural way to associate to the point $x$ a point of
the special fiber $X_0/k$ by reduction modulo $\pi$. However, by
inverting $\pi$ in the structure sheaf of $\widehat{X}$, we can
associate a generic fiber $X_\eta$ to the formal scheme
$\widehat{X}$, which is a rigid variety over $K$. Rigid geometry
provides a satisfactory theory of analytic geometry over
non-archimedean fields. A point on $X_\eta$ has coordinates in
\textit{the ring of integers} of some finite extension $K'$ of
$K$: if we denote by $R'$ the normalization of $R$ in $K'$, we can
canonically identify $X_\eta(K')$ with $\widehat{X}(R')$. By
reduction modulo $\pi$, we obtain a canonical ``contraction'' of
the generic fiber $X_\eta$ to the special fiber $X_0$. Roughly
speaking, the formal scheme $\widehat{X}$ has the advantage that
its generic and its special fiber are tightly connected, and what
glues them together are the $R'$-sections on $X$, where $R'$ runs
over the finite extensions of $R$.

A main disadvantage of rigid geometry is the artificial nature of
the topology on rigid varieties: it is not a classical topology,
but a Grothendieck topology. In the nineties, Berkovich developed
his spectral theory of non-archimedean spaces. His spaces carry a
true topology, which allows to apply classical techniques from
algebraic topology. In particular, the unit disc $R$ becomes
connected by arcs, while it is totally disconnected w.r.t. its
$\pi$-adic topology.

In Section \ref{formal}, we give a brief survey of the basic
theory of formal schemes, and Section \ref{rigid} is a crash
course on rigid geometry. Section \ref{berkovich} contains the
basic defnitions of Berkovich' approach to non-archimedean
geometry. In the final Section \ref{applic}, we briefly discuss
some applications of the theory, with special emphasis on the
relation with arc spaces of algebraic varieties, and the Milnor
fibration.

This intuitive introduction merely aims to provide some insight
into the theory of formal schemes and rigid varieties. We do not
provide proofs; instead, we chose to give a list of more thorough
introductions to the different topics dealt with in this note.
\section{Conventions and notation}
\begin{itemize}
\item For any field $F$, we denote by $F^{alg}$ an algebraic
closure, and by $F^s$ the separable closure of $F$ in $F^{alg}$.
 \item If $S$ is any scheme, a
$S$-variety is a separated reduced scheme of finite type over $S$.
\item For any locally ringed space (site) $X$, we denote the
underlying topological space (site) by $|X|$.
 \item Throughout this note, $R$ denotes a complete
discrete valuation ring, with residue field $k$, and quotient
field $K$. We fix a uniformizing parameter $\pi$, i.e. a generator
of the maximal ideal of $R$. For any integer $n\geq 0$, we denote
by $R_n$ the quotient ring $R/(\pi^{n+1})$. A finite extension
$R'$ of $R$ is by definition the normalization of $R$ in some
finite field extension $K'$ of $K$; $R'$ is again a complete
discrete valuation ring.
 \item Once we fix a value $|\pi|\in ]0,1[$, the discrete valuation $v$ on $K$ defines a
non-archimedean absolute value $|.|$ on $K$, with
$|z|=|\pi|^{v(z)}$ for $z\in K^*$. This absolute value induces a
topology on $K$, called the $\pi$-adic topology. The ideals
$\pi^nR$, $n\geq 0$, form a fundamental system of open
neighbourhoods of the zero element in $K$. The $\pi$-adic topology
is totally disconnected. The absolute value on $K$ extends
uniquely to an absolute value on $K^{alg}$. For any integer $m>0$,
we endow $(K^{alg})^m$ with the norm
$\|z\|:=\max_{i=1,\ldots,m}|z_i|$.
\end{itemize}
\section{Formal geometry}\label{formal}
In this note, we will only consider formal schemes topologically
of finite type over the complete discrete valuation ring $R$. This
case is in many respects simpler than the general one, but it
serves our purposes. For a more thorough introduction to the
theory of formal schemes, we refer to \cite[\S 10]{ega1},
\cite[no. 182]{fga}, \cite{formillusie}, or \cite{bosch}.

Intuitively, a formal scheme $X_\infty$ over $R$ consists of its
special fiber $X_0$, which is a scheme of finite type over $k$,
endowed with a structure sheaf containing additional algebraic
information on an infinitesimal neighbourhood of $X_0$.
\subsection{Affine formal schemes}
For any tuple of variables $x=(x_1,\ldots,x_m)$, we define a
 $R$-algebra $R\{x\}$ as the projective limit
$$R\{x\}:= \lim_{\stackrel{\longleftarrow}{n}}R_n[x]$$
The $R$-algebra $R\{x\}$ is canonically isomorphic to the algebra
of converging power series over $R$, i.e. the subalgebra of
$R[[x]]$ consisting of the elements
$$c(x)=\sum_{i=(i_1,\ldots,i_m)\in \N^m}\left(c_i\prod_{j=1}^{m}x_j^{i_j}\right)\ \in R[[x]]$$
 such that $c_i\to 0$ (w.r.t. the $\pi$-adic
topology on $K$) as $|i|=i_1+\ldots+i_m$ tends to $\infty$. This
means that for each $n\in\N$, there exists a value $i_0\in\N$ such
that $c_i$ is divisible by $\pi^n$ in $R$ if  $|i|>i_0$. Note that
this is exactly the condition which guarantees that the images of
$c(x)$ in the quotient rings $R_n[[x]]$ are actually polynomials,
i.e. belong to $R_n[x]$. The algebra $R\{x\}$ can also be
characterized as the sub-algebra of $R[[x]]$ consisting of those
power series which converge on the closed unit disc $R^m=\{z\in
K^m\,|\,\|z\|\leq 1\}$, since an infinite sum converges in a
non-archimedean field iff its terms tend to zero. One can show
that $R\{x\}$ is Noetherian \cite[0.(7.5.4)]{ega1}.

An $R$-algebra $A$ is called \textit{topologically of finite type}
$(tft)$ over $R$, if it is isomorphic to an algebra of the form
$R\{x_1,\ldots,x_m\}/I$, for some integer $m>0$ and some ideal
$I$. For any integer $n\geq 0$, we denote by $A_n$ the quotient
ring $A/(\pi^{n+1})$. It is a $R_n$-algebra of finite type. Then
$A$ is the limit of the projective system $(A_n)_{n\in \N}$, and
if we endow each ring $A_n$ with the discrete topology, then $A$
becomes a topological ring w.r.t. the limit topology (the
$\pi$-adic topology on $A$). By definition, the ideals $\pi^nA$,
$n>0$, form a fundamental system of open neighbourhoods of the
zero element of $A$.

To any $tft$ $R$-algebra $A$, we can associate a ringed space
$\mathrm{Spf}\,A$. It is defined as the direct limit
$$\mathrm{Spf}\,A:=\lim_{\stackrel{\longrightarrow}{n}}\mathrm{Spec}\,A_n$$
in the category of topologically ringed spaces (where the topology
on $\mathcal{O}_{\Spec A_n}$ is discrete for every $n$), so the
structure sheaf $\mathcal{O}_{\Spf A}$ is a sheaf of topological
$R$-algebras in a natural way. Moreover, one can show that the
stalks of this structure sheaf are local rings. A $tft$
\textit{affine formal $R$-scheme} is a locally topologically
ringed space in $R$-algebras which is isomorphic to a space of the
form $\Spf A$.

Note that the transition morphisms $\mathrm{Spec}\,A_m\rightarrow
\mathrm{Spec}\, A_n$, $m\leq n$, are nilpotent immersions and
therefore homeomorphisms. Hence, the underlying topological space
$|\mathrm{Spf}\,A|$ of $\Spf A$ is the set of \textit{open} prime
ideals $J$ of $A$ (i.e. prime ideals containing $\pi$), endowed
with the Zariski topology, and it is canonically homeomorphic to
$|\mathrm{Spec}\,A_0|$.

So we see that $\Spf A$ is the locally topologically ringed space
in $R$-algebras
$$(|\Spec A_0|,\lim_{\stackrel{\longleftarrow}{n}}\mathcal{O}_{\Spec A_n})$$
In particular, we have $\mathcal{O}_{\Spf A}(\Spf A)=A$. Whenever
$f$ is an element of $A$, we denote by $D(f)$ the set of open
prime ideals of $A$ which do not contain $f$. This is an open
subset of $|Spf\,A|$, and the ring of sections $\mathcal{O}_{\Spf
A}(D(f))$ is the $\pi$-adic completion $A_{\{f\}}$ of the
localization $A_f$.

 A
morphism between $tft$ affine formal $R$-schemes is by definition
a morphism of locally ringed spaces in $R$-algebras\footnote{Such
a morphism is automatically continuous w.r.t. the topology on the
structure sheaves, since it maps $\pi$ to itself; because the
topology is the $\pi$-adic one, it is determined by the
$R$-algebra structure. This is specific to so-called $R$-adic
formal schemes and does not hold for more general formal
$R$-schemes.}. If $h:A\rightarrow B$ is a morphism of $tft$
$R$-algebras, then $h$ induces a direct system of morphisms of
$R$-schemes $\Spec B_n\rightarrow \Spec A_n$ and by passage to the
limit a morphism of $tft$ affine formal $R$-schemes $Spf(h):\Spf
B\rightarrow \Spf A$. The resulting functor $\mathrm{Spf}$ induces
an equivalence between the opposite category of $tft$
$R$-algebras, and the category of $tft$ affine formal $R$-schemes,
just like in the algebraic scheme case.

The \textit{special fiber} $X_0$ of the affine formal $R$-scheme
$X_\infty=\mathrm{Spf}\,A$ is the $k$-scheme
$X_0=\mathrm{Spec}\,A_0$. As we've seen, the natural morphism of
topologically locally ringed spaces $X_0\rightarrow X_\infty$ is a
homeomorphism.

\begin{example}
Any finite extension $R'$ of $R$ is a $tft$ $R$-algebra. The
affine formal scheme $\Spf R'$ consists of a single point,
corresponding to the maximal ideal of $R'$, but the ring of
sections on this point is the entire ring $R'$. So in some sense
the infinitesimal information in the topology of $\Spec R'$ (the
generic point) is transferred to the structure sheaf of $\Spf R'$.

If $A=R\{x,y\}/(\pi-xy)$ and $X_\infty=\Spf A$, then as a
topological space $X_\infty$ coincides with its special fiber
$X_0=\Spec k[x,y]/(xy)$, but the structure sheaf of $X_\infty$ is
much ``thicker'' than the one of $X_0$. The formal $R$-scheme
$X_\infty$ should be seen as an infinitesimal tubular
neighbourhood around $X_0$.
\end{example}
\subsection{Formal schemes}
A \textit{formal scheme $X_\infty$ topologically of finite type
($tft$) over $R$} is a locally topologically ringed space in
$R$-algebras, which has a finite open cover by $tft$ affine formal
$R$-schemes.
 A
morphism between $tft$ formal $R$-schemes is a morphism of locally
ringed spaces in $R$-algebras.

It is often convenient to describe $X_\infty$ in terms of the
direct system $(X_n:=X_\infty\times_R R_n)_{n\geq 0}$. The locally
ringed space $X_n$ is a scheme of finite type over $R_n$, for any
$n$; if $X_\infty=\Spf A$, then $X_n=\Spec A_n$. For any pair of
integers $0\leq m\leq n$, the natural map of $R_n$-schemes
$u_{m,n}:X_m\rightarrow X_n$ induces an isomorphism of
$R_m$-schemes $X_m \cong X_n\times_{R_n} R_m$. The scheme $X_0$ is
called the special fiber of $X_\infty$, and $X_n$ is the $n$-th
infinitesimal neighborhood of $X_0$ in $X_\infty$ (or
``thickening''). The natural morphism of locally topologically
ringed spaces $X_n\rightarrow X_\infty$ is a homeomorphism for
each $n\geq 0$.

Conversely, if $(X_n)_{n\geq 0}$ is a direct system of $R$-schemes
of finite type such that $\pi^{n+1}=0$ on $X_n$ and such that the
transition morphism $u_{m,n}:X_m\rightarrow X_n$ induces an
isomorphism of $R_m$-schemes $X_m \cong X_n\times_{R_n} R_m$ for
each $0\leq m\leq n$, then this direct system determines a $tft$
formal $R$-scheme $X_\infty$ by putting
$$X_\infty:=\lim_{\stackrel{\longrightarrow}{n}}X_n$$ as a locally topologically
ringed space in $R$-algebras.

In the same way, giving a morphism $f:X_\infty\rightarrow
Y_\infty$ between $tft$ formal $R$-schemes amounts to giving a
series of morphisms $(f_n:X_n\rightarrow Y_n)_{n\geq 0}$, where
$f_n$ is a morphism of $R_n$-schemes, and all the squares

$$\begin{CD}
X_m@>>> X_n \\@Vf_mVV @VVf_nV \\Y_m@>>> Y_n
\end{CD}$$

\noindent commute. In other words, a morphism of $tft$ formal
$R$-schemes consists of a compatible system of morphisms between
all the infinitesimal neighbourhoods of the special fibers.

The formal scheme $X_\infty$ is called separated if the scheme
$X_n$ is separated for each $n$. In fact, this will be the case as
soon as the special fiber $X_0$ is separated. We will work in the
category of separated formal schemes, topologically of finite type
over $R$; we'll call these objects \textit{$stft$ formal
$R$-schemes}.

 A $stft$
formal scheme $X_\infty$ over $R$ is flat if its structure sheaf
has no $\pi$-torsion. A typical example of a non-flat $stft$
formal $R$-scheme is one with an irreducible component
concentrated in the special fiber. A flat $stft$ formal $R$-scheme
can be thought of as a continuous family of schemes over the
inifintesimal disc $\Spf R$.
Any $stft$ formal $R$-scheme has a maximal flat closed formal
subscheme, obtained by killing $\pi$-torsion.
\subsection{Coherent modules}
Let $A$ be a $tft$ $R$-algebra. An $A$-module $N$ is coherent iff
it is finitely generated. Any such module $N$ defines a sheaf of
modules on $\mathrm{Spf}\,A$ in the usual way. A coherent sheaf of
modules $\mathcal{N}$ on a $stft$ formal $R$-scheme $X_\infty$ is
obtained by gluing coherent modules on affine open formal
subschemes.

A more convenient description is the following: the category of
coherent sheaves $\mathcal{N}$ on $X_\infty$ is equivalent to the
category of direct systems $(\mathcal{N}_n)_{n\geq 0}$, where
$\mathcal{N}_n$ is a coherent sheaf on the scheme $X_n$, and the
$\mathcal{O}_{X_n}$-linear transition map
$v_{m,n}:\mathcal{N}_m\rightarrow \mathcal{N}_n$ induces an
isomorphism of coherent $\mathcal{O}_{X_m}$-modules $\mathcal{N}_m
\cong u_{m,n}^*\mathcal{N}_n$ for any pair $m\leq n$. Morphisms
between such systems are defined in the obvious way.

\subsection{The completion functor}\label{completion}
Let $X$ be any Noetherian scheme and $\mathcal{J}$ a coherent
ideal sheaf on $X$, and denote by $V(\mathcal{J})$ the closed
subscheme of $X$ defined by $\mathcal{J}$. The $\mathcal{J}$-adic
completion $\widehat{X/\mathcal{J}}$of $X$ is the limit of the
direct system of schemes $(V(\mathcal{J}^n))_{n>0}$ in the
category of topologically ringed spaces (where
$\mathcal{O}_{V(\mathcal{J}^n)}$ carries the discrete topology).
This is, in fact, a formal scheme, but in general not of the kind
we have defined before; we include the construction here for later
use. If $h:Y\rightarrow X$ is a morphism of Noetherian schemes,
and if we denote by $\mathcal{K}$ the inverse image
$\mathcal{J}\mathcal{O}_{Y}$ of $\mathcal{J}$ on $Y$, then $h$
defines a direct system of morphisms of schemes
$V(\mathcal{K}^n)\rightarrow V(\mathcal{J}^n)$ and by passage to
the limit a morphism of topologically locally ringed spaces
$\widehat{Y/\mathcal{K}}\rightarrow \widehat{X/\mathcal{J}}$,
called the $\mathcal{J}$-adic completion of $h$.

If $X$ is a separated $R$-scheme of finite type and $\mathcal{J}$
is the ideal generated by $\pi$, then the $\mathcal{J}$-adic
completion of $X$ is the limit of the direct system
$(X_n=X\times_R R_n)_{n\geq 0}$, and this is a $stft$ formal
$R$-scheme which we denote simply by $\widehat{X}$. It is called
the \textit{formal ($\pi$-adic) completion} of the $R$-scheme $X$.
Its special fiber $X_0$ is canonically isomorphic to the fiber of
$X$ over the closed point of $\mathrm{Spec}\,R$. The formal scheme
$\widehat{X}$ is flat iff $X$ is flat over $R$. Intuitively,
$\widehat{X}$ should be seen as the infinitesimal tubular
neighbourhood of $X_0$ in $X$. As a topological space, it
coincides with $X_0$, but additional infinitesimal information is
contained in the structure sheaf.

\begin{example}
If $X=\mathrm{Spec}\,R[x_1,\ldots,x_n]/(f_1,\ldots,f_\ell)$, then
its formal completion $\widehat{X}$ is simply
$\mathrm{Spf}\,R\{x_1,\ldots,x_n\}/(f_1,\ldots,f_\ell)$.
\end{example}

By the above construction, a morphism of separated $R$-schemes of
finite type $f:X\rightarrow Y$ induces a morphism of formal
$R$-schemes $\hat{f}:\widehat{X}\rightarrow \widehat{Y}$ between
the formal $\pi$-adic completions of $X$ and $Y$. We get a
completion functor
$$\widehat{ }\,:(sft-Sch/R)\rightarrow (stft-For/R):X\mapsto \widehat{X}$$
where $(sft-Sch/R)$ denotes the category of separated $R$-schemes
of finite type, and $(stft-For/R)$ denotes the category of
separated formal schemes, topologically of finite type over $R$.

For a general pair of separated $R$-schemes of finite type $X,Y$,
the completion map
$$C_{X,Y}:Hom_{(sft-Sch/R)}(X,Y)\rightarrow Hom_{(stft-For/R)}(\widehat{X},\widehat{Y}):f\mapsto\widehat{f}$$
is injective, but not bijective. It is a bijection, however, if
$X$ is proper over $R$: this is a corollary of Grothendieck's
Existence Theorem; see \cite[5.4.1]{ega3}. In particular, the
completion map induces a bijection between $R'$-sections of $X$,
and $R'$-sections of $\widehat{X}$ (i.e. morphisms of formal
$R$-schemes $\Spf R'\rightarrow \widehat{X}$), for any finite
extension $R'$ of the complete discrete valuation ring $R$.
Indeed: $\Spec R'$ is a finite, hence proper $R$-scheme, and its
formal $\pi$-adic completion is $\Spf R'$.

\begin{example}
If $X=\Spec B$, with $B$ an $R$-algebra of finite type, and
$Y=\Spec R[z]$, then
$$Hom_{(sft-Sch/R)}(X,Y)=B$$
On the other hand, if we denote by $\widehat{B}$ the $\pi$-adic
completion of $B$, then $\widehat{X}=\Spf \widehat{B}$ and
$\widehat{Y}=\Spf R\{z\}$, and we find
$$Hom_{(stft-For/R)}(\widehat{X},\widehat{Y})=\widehat{B}$$
The completion map $C_{X,Y}$ is given by the natural injection
$B\rightarrow \widehat{B}$; it is not surjective in general, but
it is surjective if $B$ is finite over $R$.
\end{example}

If $X$ is a separated $R$-scheme of finite type, and $\mathcal{N}$
is a coherent sheaf of $\mathcal{O}_X$-modules, then $\mathcal{N}$
induces a direct system $(\mathcal{N}_n)_{n\geq 0}$, where
$\mathcal{N}_n$ is the pull-back of $\mathcal{N}$ to $X_n$. This
system defines a coherent sheaf of modules $\widehat{\mathcal{N}}$
on $\widehat{X}$. If $X$ is proper over $R$, it follows from
Grothendieck's Existence Theorem that the functor
$\mathcal{N}\rightarrow \widehat{\mathcal{N}}$ is an equivalence
between the category of coherent $\mathcal{O}_{X}$-modules and the
category of coherent $\mathcal{O}_{\widehat{X}}$-modules
\cite[5.1.6]{ega3}. Moreover, there is a canonical isomorphism
$H^q(\widehat{X},\widehat{\mathcal{N}})\cong H^q(X,\mathcal{N})$
for each coherent $\mathcal{O}_{X}$-module $\mathcal{N}$ and each
$q\geq 0$.

If a $stft$ formal $R$-scheme $Y_\infty$ is isomorphic to the
$\pi$-adic completion $\widehat{Y}$ of a separated $R$-scheme $Y$
of finite type, we call the formal scheme $Y_\infty$
\textit{algebrizable}, with algebraic model $Y$. The following
theorem is the main criterion to recognize algebrizable formal
schemes \cite[5.4.5]{ega3}: if $Y_\infty$ is proper over $R$, and
$\mathcal{L}$ is an invertible $\mathcal{O}_{Y_\infty}$-bundle
such that the pull-back $\mathcal{L}_0$ of $\mathcal{L}$ to $Y_0$
is ample, then $Y_\infty$ is algebrizable. Moreover, the algebraic
model $Y$ for $Y_\infty$ is unique up to canonical isomorphism,
there exists a unique line bundle $\mathcal{M}$ on $Y$ with
$\mathcal{L}=\widehat{\mathcal{M}}$, and $\mathcal{M}$ is ample.
For an example of a proper formal $\C[[t]]$-scheme which is not
algebrizable, see \cite[5.24(b)]{formillusie}.
\subsection{Formal blow-ups}
 Let $X_\infty$ be a flat $stft$ formal $R$-scheme, and let
 $\mathcal{I}$ be a coherent ideal sheaf on $X_\infty$ such that $\mathcal{I}$
 contains a power of the uniformizing parameter $\pi$.
 We can define the formal blow-up of $X_\infty$ at the center $\mathcal{I}$ as follows
 \cite[\S\,2]{formrigI}: if $X_\infty=\Spf A$ is affine, and $I$
 is the ideal of global sections of $\mathcal{I}$ on $X_\infty$,
 then the formal blow-up of $X_\infty$ at $\mathcal{I}$ is the
 $\pi\mathcal{O}_{\Spec A}$-adic completion of the blow-up of $\Spec A$ at $I$. The
 general case is obtained by gluing.

  The formal blow-up of $X_\infty$ at $\mathcal{I}$ is again a flat $stft$ formal
  $R$-scheme,and
the composition of two formal blow-ups is again a
 formal blow-up \cite[2.1+2.5]{formrigI}.
If $X$ is a separated $R$-scheme of finite type and $\mathcal{I}$
is a coherent ideal sheaf on $X$ containing a power of $\pi$, then
the formal blow-up of $\widehat{X}$ at $\widehat{\mathcal{I}}$ is
canonically isomorphic to the $\pi$-adic completion of the blow-up
of $X$ at $\mathcal{I}$.

\section{Rigid geometry}\label{rigid}
In this note, we'll only be able to cover the basics of rigid
geometry. We refer the reader to the books \cite{BGR,fresnel} and
the research papers \cite{bert,formrigI,Raynaud,tate} for a more
thorough introduction. A nice survey on Tate's approach to rigid
geometry can be found in \cite{liu-rigid}.

\subsection{Analytic geometry over non-Archimedean fields}
Let $L$ be a non-archimedean field (i.e. a field which is complete
w.r.t. an absolute value which satisfies the ultrametric
property); we assume that the absolute value on $L$ is
non-trivial. For instance, if $K$ is our complete discretely
valued field, then we can turn $K$ into a non-archimedean field by
fixing a value $|\pi|\in ]0,1[$ and putting $|x|=|\pi|^{v(x)}$ for
$x\in K^*$, where $v$ denotes the discrete valuation on $K$ (by
convention, $v(0)=\infty$ and $|0|=0$).

The absolute value on $L$ extends uniquely to any finite extension
of $L$, and hence to $L^s$ and $L^{alg}$. We denote by
 $\widehat{L^{alg}}$ the completion of $L^{alg}$, and by $\widehat{L^s}$ the
 closure of $L^s$ in $\widehat{L^{alg}}$; these
are again non-archimedean fields. We denote by $L^o$ the valuation
ring $\{x\in L\,|\,|x|\leq 1\}$, by $L^{oo}$ its maximal ideal
$\{x\in L\,|\,|x|<1\}$, and by $\widetilde{L}$ the residue field
$L^o/L^{oo}$. For $L=K$ we have $L^o=R$, $L^{oo}=(\pi)$ and
$\widetilde{L}=k$.


Since $L$ is endowed with an absolute value, one can use this
structure to develop a theory of analytic varieties over $L$ by
mimicing the construction over $\C$. Na\"ively, we can define
analytic functions on open subsets of $L^n$ as $L$-valued
functions which are locally
 defined by a converging
power series with coefficients in $L$. However, we are immediately
confronted with some pathological phenomena. Consider, for
instance, the $p$-adic unit disc
$$\Z_p=\{x\in \Q_p\,|\,|x|\leq 1\}$$ The partition
$$\{p\,\Z_p, 1+p\,\Z_p,\ldots,(p-1)+p\,\Z_p\}$$
is an open cover of $\Z_p$ w.r.t. the $p$-adic topology. Hence,
the characteristic function of $p\,\Z_p$ is analytic, according to
our na\"ive definition. This contradicts some elementary
properties that one expects an analytic function to have. The
cause of this and similar pathologies, is the fact that the unit
disc $\Z_p$ is totally disconnected with respect to the $p$-adic
topology. In this approach, there are ``too many'' analytic
functions, and ``too few'' analytic varieties (for instance, with
this definition, any compact $p$-adic manifold is isomorphic to a
disjoint union of $i$ unit discs, with $i\in \{0,\ldots,p-1\}$ its
Serre invariant \cite{serre}).

Rigid geometry is as a more refined approach to non-Archimedean
analytic geometry, turning the unit disc into a connected space.
Rigid spaces are endowed with a certain Grothendieck topology,
only allowing a special type of covers.

We'll indicate two possible approaches to the theory of rigid
varieties over $L$. The first one is due to Tate \cite{tate}, the
second one to Raynaud \cite{Raynaud}. If we return to our example
of the $p$-adic unit disc $\Z_p$, Tate's construction can be
understood as follows. In fact, we already know what the
``correct'' algebra of analytic functions on $\Z_p$ should be: the
 power series with coefficients in $K$ which converge
 \textit{globally} on $\Z_p$. Tate's idea is to start from this
 algebra, and then to construct a space on which these functions
 naturally live. This is similar to the construction of the spectrum
 of a ring in algebraic geometry. Raynaud observed that a certain class of
  Tate's
 rigid varieties can be characterized in terms of formal
 schemes.
\subsection{Tate algebras}\label{subsec-tate}
 The basic objects in Tate's theory are the \textit{algebras of
 converging power series over $L$}
$$T_m=L\{x_1,\ldots,x_m\}=\{\alpha=\sum_{i\in \N^m}(\alpha_i \prod_{j=1}^{m}x_j^{i_j})\in L[[x_1,\ldots,x_m]]\,|\,
|\alpha_i|\rightarrow 0\ \mbox{as}\ |i|\rightarrow \infty\} $$
where $|i|=\sum_{j=1}^mi_j$. The convergence condition implies in
particular that for each $\alpha$ there exists $i_0\in\N$ such
that for $|i|>i_0$ the coefficient $\alpha_i$ belongs to $L^o$.
The algebra $T_m$ is the algebra of power series over $L$ which
converge on the closed unit polydisc $(L^o)^m$ in $L^m$ (since an
infinite sum converges in a non-archimedean field iff its terms
tend to zero). Note that, for $L=K$, $T_m\cong
R\{x_1,\ldots,x_m\}\otimes_R K$. Analogously, we can define an
algebra of converging power series $B\{x_1,\ldots,x_m\}$ for any
Banach algebra $B$. The algebra $T_m$ is a Banach algebra for the
sup-norm $\|f\|_{sup}=\max_i |\alpha_i|$. It is Noetherian, and
any ideal $I$ is closed, so that the quotient $T_m/I$ is again a
Banach algebra w.r.t. the residue norm.

A Tate algebra, or $L$-affinoid algebra, is a $L$-algebra $A$
isomorphic to such a quotient $T_m/I$. The residue norm on $A$
depends on the presentation $A\cong T_m/I$. However, any morphism
of $L$-algebras $T_m/I\rightarrow T_m/J$ is automatically
continuous, so in particular, the residue norm on $A$ is
well-defined up to equivalence, and the induced topology on $A$ is
independent of the chosen presentation. For any maximal ideal $y$
of $A$, the residue field $A/y$ is a finite extension of $L$. For
proofs of all these facts, we refer to \cite[3.2.1]{fresnel}.

By Proposition 1 of \cite[7.1.1]{BGR}, the maximal ideals $y$ of
$T_m$ correspond bijectively to $G(L^{alg}/L)$-orbits of tuples
$(z_1,\ldots,z_m)$, with $z_i\in (L^{alg})^o$, via the map
$$y\mapsto \{(\varphi(x_1),\ldots,\varphi(x_m))\,|\,\varphi:T_m/y\hookrightarrow L^{alg}\}$$
where $\varphi$ runs through the $L$-embeddings of $T_m/y$ in
$L^{alg}$. In particular, for any morphism of $L$-algebras
$\psi:T_m\rightarrow L^{alg}$ and any index $i$, the element
$\psi(x_i)$ belongs to $(L^{alg})^o$. It follows that $\psi$ is
contractive, in the sense that $|\psi(a)|\leq \|a\|_{sup}$ for any
$a$ in $T_m$.

The fact that we obtain tuples of elements in $(L^{alg})^o$,
rather than $L^{alg}$, might look strange at first; it is one of
the most characteristic properties of Tate's rigid varieties. Let
us consider an elementary example. If $z$ is an element of $L$,
then $x-z$ is invertible in $L\{x\}$ iff $z$ does not belong to
$L^o$. Indeed, for $z\neq 0$ the coefficients of the formal power
series $1/(x-z)=-(1/z)\sum_{i\geq 0}(x/z)^i$ tend zo zero iff
$|z|>1$, i.e. iff $z\notin L^o$. So $(x-z)$ only defines a maximal
ideal in $L\{x\}$ if $z\in L^o$.
\subsection{Affinoid spaces}\label{subsec-aff}
The category of \textit{$L$-affinoid spaces} is by definition the
opposite category of the category of Tate algebras over $L$.  For
any $L$-affinoid space $X$, we'll denote the corresponding Tate
algebra by $A(X)$, and we call it the \textit{algebra of analytic
functions} on $X$. Conversely, for any Tate algebra $A$, we denote
the corresponding affinoid space by $\mathrm{Sp}\,A$ (some authors
use the notation $\mathrm{Spm}$ instead). For any $m\geq 0$, the
affinoid space $\Sp T_m$ is called the closed unit disc of
dimension $m$ over $L$.

To any $L$-affinoid space $X=\mathrm{Sp}\,A$, we associate the set
$X^{\flat}$ of maximal ideals of the Tate algebra $A=A(X)$. If we
present $A$ as a quotient $T_m/(f_1,\ldots,f_n)$, then elements of
$(\mathrm{Sp}\,A)^{\flat}$ correspond bijectively to
$G(L^{alg}/L)$-orbits of tuples $z=(z_1,\ldots,z_m)$, with $z_i\in
(L^{alg})^o$, and $f_j(z)=0$ for each $j$. In particular, if $L$
is algebraically closed and $X$ is the closed unit disc $\Sp T_1$,
then
$$X^{\flat}=L^o=\{x\in L\,|\,|x|\leq 1\}$$

We've seen above that, for any maximal ideal $x$ of $A$, the
quotient $A/x$ is a finite extension of $L$, so it carries a
unique prolongation of the absolute value $|.|$ on $L$. Hence, for
any $f\in A$ and any $x\in (\Sp A)^{\flat}$, we can speak of the
value $f(x)$ of $f$ at $x$ (the image of $f$ in $A/x$), and its
absolute value $|f(x)|$. In this way, elements of $A$ are viewed
as functions on $(\Sp A)^{\flat}$. Note that, if $x$ is a prime
ideal of $A$, there is in general no canonical way to extend the
absolute value on $L$ to the extension $A/x$. This is one of the
reasons for working with the maximal spetrum $(\Sp A)^{\flat}$,
rather than the prime spectrum $\Spec A$. In Berkovich' theory
(Section \ref{berkovich}) the notion of point is generalized by
admitting any prime ideal $x$ \textit{and} specifying an extension
of the absolute value on $L$ to $A/x$.

The \textit{spectral semi-norm} on $A$ is defined by
$$\|f\|_{sup}:=\sup_{x\in X^{\flat}}\,|f(x)|$$
It is a norm iff $A$ is reduced. By the \textit{maximum modulus
principle} \cite[6.2.1.4]{BGR}, this supremum is, in fact, a
maximum, i.e. there is a point $x$ in $X^{\flat}$ with
$|f(x)|=\|f\|_{sup}$. Moreover, for $A=T_m$, this definition
concides with the one in the previous section, by
\cite[5.1.4.6]{BGR}.

If $\varphi:A\rightarrow B$ is a morphism of $L$-affinoid
algebras, then for any maximal ideal $x$ in $B$, $\varphi^{-1}(x)$
is a maximal ideal in $A$, since $B/x$ is a finite extension of
$L$. Hence, any morphism of $L$-affinoid spaces $h:X\rightarrow Y$
induces a map $h^{\flat}:X^{\flat}\rightarrow Y^{\flat}$ on the
associated sets.
A morphism of $L$-affinoid
spaces $h:X\rightarrow Y$ is called a \textit{closed immersion} if
the corresponding morphism of $L$-affinoid algebras
$A(Y)\rightarrow A(X)$ is surjective.

We could try to endow $X^{\flat}$ with the inital topology w.r.t.
the functions $x\mapsto |f(x)|$, where $f$ varies in $A$. If $L$
is algebraically closed and if we identify $(\Sp L\{x\})^{\flat}$
with $L^o$, then this topology is simply the topology on $L^o$
defined by the absolute value. It is totally disconnected, so it
does not have the nice properties we are looking for.
\subsection{Open covers}
A morphism $h:Y\rightarrow X$ of $L$-affinoid spaces is called an
\textit{open immersion}, if it satisfies the following universal
property: for any morphism $g:Z\rightarrow X$ of $L$-affinoid
spaces such that the image of $g^{\flat}$ is contained in the
image of $h^{\flat}$ in $X^{\flat}$, there is a unique morphism
$g':Z\rightarrow Y$ such that $g=h\circ g'$. If $h$ is an open
immersion, the image $D$ of $h^{\flat}$ in $X^{\flat}$ is called
an \textit{affinoid domain}. One can show that the map $h^{\flat}$
is always injective \cite[7.2.2.1]{BGR}, so it identifies the set
$D$ with $Y^{\flat}$. The $L$-affinoid space $Y$ and the open
immersion $h:Y\rightarrow X$ are uniquely determined by the
affinoid domain $D$, up to canonical isomorphism. With slight
abuse of notation, we will identify the affinoid domain $D$ with
the $L$-affinoid space $Y$, so that we can think of an affinoid
domain as an affinoid space sitting inside $X$, and we can speak
of the Tate algebra $A(D)$ of analytic functions on $D$. If $E$ is
a subset of $D$, then $E$ is an affinoid domain in $D$ iff it is
an affinoid domain in $X$. In this case, the universal property
yields a restriction map $A(D)\rightarrow A(E)$. The intersection
of two affinoid domains is again an affinoid domain, but this does
not always hold for their union. If $h:Z\rightarrow X$ is a
morphism of $L$-affinoid spaces, then the inverse image of an
affinoid domain in $X$ is an affinoid domain in $Z$ .

\begin{example}\label{ex-affinoid}
Consider the closed unit disc $X=\Sp L\{x\}$. For $a\in L^o$ and
$r$ in the value group $|L^*|$, we denote by $D(a,r)$ the ``closed
disc'' $\{z\in X^{\flat}\,|\,|x(z)-a|\leq r\}$, and by
$D^{-}(a,r)$ the ``open disc'' $\{z\in X^{\flat}\,|\,|x(z)-a|<
r\}$. We will see below that the disc $D(a,r)$ is an affinoid
domain in $X$, with $A(D(a,r))= L\{x,T\}/(x-a-\rho T)$, where
$\rho$ is any element of $L$ with $|\rho|=r$. On the other hand,
the disc $D^{-}(a,r)$ can not be an affinoid domain in $X$, since
the function $|x(.)-a|$ does not reach its maximum on
$D^{-}(a,r)$.

Assume now that $L$ is algebraically closed. By Theorem $2$ in
\cite[9.7.2]{BGR}
 the affinoid domains in $X$
 are the finite disjoint unions of
subsets of the form
$$D(a_0,r_0)\setminus \cup_{i=1}^q D^{-}(a_i,r_i)$$
with $a_i$ in $L^o$ and $r_i$ in $|L^*|\cap ]0,1]$ for
$i=0,\ldots,q$.
\end{example}

An \textit{affinoid cover} of $X$ is a \textit{finite} set of open
immersions $u_i:U_i\rightarrow X$ such that the images of the maps
$(u_i)^{\flat}$ cover $X^{\flat}$. A special kind of affinoid
covers is constructed as follows: take analytic functions
$f_1,\ldots,f_n$ in $A(X)$, and suppose that these elements
generate the unit ideal $A(X)$. Consider, for each $i=1,\ldots,n$,
the $L$-affinoid space $U_i$ given by
$$A(U_i)=A(X)\{T_1,\ldots,T_n\}/(f_j-T_jf_i)_{j=1,\ldots,n}$$
The obvious morphism $u_i:U_i\rightarrow X$ is an open immersion,
and $U_i$ is called a \textit{rational subspace} of $X$. The image
of $(u_i)^{\flat}$ is the set of points $x$ of $X^{\flat}$ such
that $|f_i(x)|\geq |f_j(x)|$ for $j=1,\ldots,n$. Indeed: using the
fact that a morphism of $L$-algebras $T_m\rightarrow L^{alg}$ is
contractive (Section \ref{subsec-tate}), and the assumption that
$f_1,\ldots,f_n$ generate $A(X)$, one shows that a morphism of
$L$-algebras $\psi:A(X)\rightarrow L^{alg}$ factors through a
morphism of $L$-algebras $\psi_i:A(U_i)\rightarrow L^{alg}$ iff
$\psi(f_i)\neq 0$ and $\psi(T_j)=\psi(f_j)/\psi(f_i)$ belongs to
$(L^{alg})^o$, i.e. $|\psi(f_j)|\leq |\psi(f_i)|$. In this case,
$\psi_i$ is unique.

The set of morphisms $\{u_1,\ldots,u_n\}$ is an affinoid cover,
and is called a \textit{standard cover}. It is a deep result that
any affinoid domain of $X$ is a finite union of rational subsets
of $X$, and any affinoid cover of $X$ can be refined by a standard
cover \cite[7.3.5.3+8.2.2.2]{BGR}.

One of the cornerstones in the theory of rigid varieties is Tate's
Acyclicity Theorem \cite[8.2.1.1]{BGR}. It states that analytic
functions on any affinoid cover $\{u_i:U_i\rightarrow X\}_{i\in
I}$ satisfy the gluing property: the sequence
$$A(X)\rightarrow \prod_{i\in I}A(U_i)\rightrightarrows
\prod_{(i,j)\in I^2}A(U_i\cap U_j)$$ is exact.

Now we can define, for each $L$-affinoid space $X$, a topology on
the associated set $X^{\flat}$. It will not be a topology in the
classical sense, but a Grothendieck topology, a generalization of
the topological concept in the framework of categories. A
Grothendieck topology specifies a class of opens (admissible
opens) and, for each admissible open, a class of covers
(admissible covers). These have to satsify certain axioms which
allow to develop a theory of sheaves and cohomology in this
setting. A space with a Grothendieck topology is called a
\textit{site}. Any topological space (in the classical sense) can
be viewed as a site in a canonical way: the admissible opens and
admissible covers are the open subsets and the open covers. For
our purposes, we do not need the notion of Grothendieck topology
in its most abstract and general form: a sufficient treatment is
given in \cite[9.1.1]{BGR}.

The \textit{weak} $G$-topology on an $L$-affinoid space $X$ is
defined as follows: the admissible open sets of $X^{\flat}$ are
the affinoid domains, and the admissible covers are the affinoid
covers \cite[9.1.4]{BGR}. Any morphism $h$ of $L$-affinoid spaces
is continuous w.r.t. the weak $G$-topology (meaning that the
inverse image under $h^{\flat}$ of an admissible open is again an
admissible open, and the inverse image of an admissible cover is
again an admissible cover). We can define a presheaf of
$L$-algebras $\mathcal{O}_X$ on $X^{\flat}$ with respect to this
topology, by putting $\mathcal{O}_X(D)=A(D)$ for any affinoid
domain $D$ of $X$ (with the natural restriction maps). By Tate's
Acyclicity Theorem, $\mathcal{O}_X$ is a sheaf. Note that the
exact definition of the weak $G$-topology varies in literature:
sometimes the admissible opens are taken to be the finite unions
of rational subsets in $X$, and the admissible covers are the
covers by admissible opens with a finite subcover (e.g. in
\cite[\S\,4.2]{fresnel}).

In the theory of Grothendieck topologies, there is a canonical way
to refine the topology without changing the associated category of
sheaves \cite[9.1.2]{BGR}. This refinement is important to get
good gluing properties for affinoid spaces, and to obtain
continuity of the analytification map (Section \ref{subsec-an}).
This leads to the following definition of the \textit{strong}
$G$-topology on a $L$-affinoid space $X$. \begin{itemize} \item
The admissible open sets are (possibly infinite) unions
$\cup_{i\in I} D_i$ of affinoid domains $D_i$ in $X$, such that,
for any morphism of $L$-affinoid spaces $h:Y\rightarrow X$, the
image of $h^{\flat}$ in $X^{\flat}$ is covered by a finite number
of $D_i$. \item An admissible cover of an admissible open subset
$V\subset X^{\flat}$ is a (possibly infinite) set of admissible
opens $\{V_j\,|\,j\in J\}$ in $X^{\flat}$ such that $V=\cup_j
V_j$, and such that, for any morphism of $L$-affinoid spaces
$\varphi:Y\rightarrow X$ with $Im\,(\varphi^{\flat})\subset V$,
the cover $\{(\varphi^{\flat})^{-1}(V_j)\,|\,j\in J\}$ of $Y$ can
be refined by an affinoid cover.
\end{itemize}

Any morphism of $L$-affinoid spaces is continuous w.r.t. the
strong $G$-topology. The strong $G$-topology on $X=\Sp A$ is finer
than the Zariski topology on the maximal spectrum of $A$ (this
does not hold for the weak $G$-topology). From now on, we'll endow
all $L$-affinoid spaces $X$ with the strong $G$-topology. The
strucure sheaf $\mathcal{O}_X$ of $X$ extends uniquely to a sheaf
of $L$-algebras w.r.t. the strong $G$-topology, which is called
the sheaf of analytic functions on $X$. One can show that its
stalks are local rings. In this way, we associate to any
$L$-affinoid space $X$ a locally ringed site in $L$-algebras
$(X^{\flat},\mathcal{O}_{X})$.

For any morphism of $L$-affinoid spaces $h:Y\rightarrow X$, there
is a morphism of sheaves of $L$-algebras $\mathcal{O}_X\rightarrow
(h^{\flat})_*\mathcal{O}_Y$ which defines a morphism of locally
ringed spaces $(Y^{\flat},\mathcal{O}_{Y})\rightarrow
(X^{\flat},\mathcal{O}_{X})$ (if $D=\Sp A$ is an affinoid domain
in $X$, then $(h^{\flat})^{-1}(D)$ is an affinoid domain $\Sp B$
in $Y$ and there is a natural morphism of $L$-algebras
$A\rightarrow B$ by the universal property defining affinoid
domains). This construction defines a functor from the category of
$L$-affinoid spaces to the category of locally ringed spaces in
$L$-algebras, and this functor is fully faithful
\cite[9.3.1.1]{BGR}, i.e. every morphism of locally ringed sites
in $L$-algebras $((\Sp B)^{\flat},\mathcal{O}_{\Sp B})\rightarrow
((\Sp A)^{\flat},\mathcal{O}_{\Sp A})$ is induced by a morphism of
$L$-algebras $A\rightarrow B$. With slight abuse of notation, we
will call the objects in its essential image also $L$-affinoid
spaces, and we'll identify an $L$-affinoid space $X$ with its
associated locally ringed site in $L$-algebras
$(X^{\flat},\mathcal{O}_{X})$.

If $D$ is an affinoid domain in $X$, then the strong $G$-topology
on $X$ restricts to the strong $G$-topology on $D$, and the
restriction of $\mathcal{O}_{X}$ to $D$ is the sheaf of analytic
functions $\mathcal{O}_{D}$.

One can check that the affinoid space $\Sp T_m$ is connected with
respect to the strong and the weak $G$-topology, for any $m\geq
0$. More generally, connectedness of an $L$-affinoid space $X=\Sp
A$ is equivalent for the weak $G$-topology, the strong
$G$-topology, and the Zariski topology
\cite[9.1.4,\,Prop.\,8]{BGR}, and it is also equivalent to the
property that the ring $A$ has no non-trivial idempotents; so the
$G$-topologies nicely reflect the algebraic structure of $A$.

\begin{example}\label{strongtop}
Let $X$ be the closed unit disc $\Sp L\{x\}$. The set
$$U_1=\{z\in X\,|\,|x(z)|=1\}=\{z\in X\,|\,|x(z)|\geq 1\}$$ is a rational domain in $X$, so it
is an admissible open already for the weak $G$-topology. The
algebra of analytic functions on $U_1$ is given by
$$\mathcal{O}_X(U_1)=L\{x,T\}/(xT-1)$$

The set
$$U_2=\{z\in X\,|\,|x(z)|<1\}$$ is not an admissible open for the
weak $G$-topology (it cannot be affinoid since the function
$|x(.)|$ does not reach a maximum on $U_2$), but it is an
admissible open for the strong $G$-topology: we can write it as an
infinite union  of rational domains
$$U_2^{(n)}=\{z\in X\,|\,|x(z)|^n\leq |a|\}$$ where $n$ runs throug $\N^\ast$ and $a$ is any
non-zero element of $L^{oo}$.

This family satisfies the finiteness condition in the definition
of the strong $G$-topology: if $Y\rightarrow X$ is any morphism of
$L$-affinoid spaces whose image is contained in $U_2$, then by the
maximum principle (Section \ref{subsec-aff}) the pull-back of the
function $|x(.)|$ to $Y$ reaches its maximum on $Y$, so the image
of $Y$ is contained in $U_2^{(n)}$ for $n$ sufficiently large.

The algebra $\mathcal{O}_X(U_2)$ of analytic functions on $U_2$
consists of the elements $\sum_{i\geq 0} a_ix^i$ of $L[[x]]$ such
that $|a_i|r^i$ tends to zero as $i\to\infty$, for any $r\in
]0,1[$.

Hence, we can write $X$ as a disjoint union $U_1\sqcup U_2$ of
admissible opens. This does not contradict the fact that $X$ is
connected, because $\{U_1,U_2\}$ is \textit{not} an admissible
cover, since it can not be refined by a (finite!) affinoid cover.
\end{example}
\subsection{Rigid varieties}
Now, we can give the definition of a general rigid variety over
$L$. It is a set $X$, endowed with a Grothendieck
topology\footnotetext{To be precise, this Grothendieck topology
should satisfy certain additional axioms; see
\cite[9.3.1.4]{BGR}.} and a sheaf of $L$-algebras $\mathcal{O}_X$,
such that $X$ has an admissible cover $\{U_i\}_{i\in I}$ with the
property that each locally ringed space
$(U_i,\mathcal{O}_{X}|_{U_i})$ is isomorphic to an $L$-affinoid
space. An admissible open $U$ in $X$ is called an affinoid domain
in $X$ if $(U,\mathcal{O}_{X}|_U)$ is isomorphic to an
$L$-affinoid space. If $X$ is affinoid, this definition is
compatible with the previous one. A morphism $Y\rightarrow X$ of
rigid varieties over $L$ is a morphism of locally ringed spaces in
$L$-algebras.

 A rigid variety over
$L$ is called \textit{quasi-compact} if it is a finite union of
affinoid domains. It is called \textit{quasi-separated} if the
intersection of any pair of affinoid domains is quasi-compact, and
\textit{separated} if the diagonal morphism is a closed immersion.
\subsection{Analytification of a $L$-variety}\label{subsec-an}
For any $L$-scheme $X$ of finite type, we can endow the set
$X^{o}$ of closed points of $X$ with the structure of a rigid
$L$-variety.

More precisely, by \cite[0.3.3]{bert} and \cite[5.3]{liu-rigid}
there exists a functor
$$(\,.\,)^{an}:(ft-Sch/L)\rightarrow (Rig/L)$$
from the category of $L$-schemes of finite type, to the category
of rigid $L$-varieties, such that
\begin{enumerate}
\item For any $L$-scheme of finite type $X$, there exists a
natural morphism of locally ringed sites
$$i:X^{an}\rightarrow X$$ which induces a bijection
between the underlying set of $X^{an}$ and the set $X^o$ of closed
points of $X$. The couple $(X^{an},i)$ satisfies the following
universal property: for any rigid variety $Z$ over $L$ and any
morphism of locally ringed sites $j:Z\rightarrow X$, there exists
a unique morphism of rigid varieties $j':Z\rightarrow X^{an}$ such
that $j=i\circ j'$.

\item If $f:X'\rightarrow X$ is a morphism of $L$-schemes of
finite type, the square
$$\begin{CD}
(X')^{an}@>f^{an}>>X^{an}
\\@ViVV @VViV
\\ X' @>>> X
\end{CD}$$
commutes.

 \item The functor $(\,.\,)^{an}$ commutes with fibered
products, and takes open (resp. closed) immersions of $L$-schemes
to open (resp. closed) immersion of rigid $L$-varieties. In
particular, $X^{an}$ is separated if $X$ is separated.
\end{enumerate}
We call $X^{an}$ the analytification of $X$. It is quasi-compact
if $X$ is proper over $L$, but not in general. The analytification
functor has the classical GAGA properties: if $X$ is proper over
$L$, then analytification induces an equivalence between coherent
$\mathcal{O}_X$-modules and coherent
$\mathcal{O}_{X^{an}}$-modules, and the cohomology groups agree; a
closed rigid subvariety of $X^{an}$ is the analytification of an
algebraic subvariety of $X$; and for any $L$-variety $Y$, all
morphisms $X^{an}\rightarrow Y^{an}$ are algebraic. These results
can be deduced from Grothendieck's Existence Theorem; see
\cite[2.8]{formrig0}.

\begin{example}\label{ex-rigline}
Let $D$ be the closed unit disc $\Sp L\{x\}$, and consider the
endomorphism $\sigma$ of $D$ mapping $x$ to $a\cdot x$, for some
non-zero $a\in L^{oo}$. Then $\sigma$ is an isomorphism from $D$
onto the affinoid domain $D(0,|a|)$ in $D$ (notation as in Example
\ref{ex-affinoid}).
 The rigid affine line
$(\A^1_L)^{an}$ is the limit of the direct system
$$\begin{CD}D@>\sigma >> D@>\sigma>> \ldots\end{CD}$$ in the category of
locally ringed sites in $L$-algebras.
 Intuitively, it is obtained as the union of an
infinite number of concentric closed discs whose radii tend to
$\infty$.
\end{example}
\subsection{Rigid spaces and formal schemes}
Finally, we come to a second approach to the theory of rigid
spaces, due to Raynaud \cite{Raynaud}. We will only deal with the
case where $L=K$ is a complete discretely valued field, but the
theory is valid in greater generality (see \cite{formrigI}).

We've seen before that the underlying topological space of a
$stft$ formal $R$-scheme $X_\infty$ coincides with the underlying
space of its special fiber $X_0$. Nevertheless, the structure
sheaf of $X_\infty$ contains information on an infinitesimal
neighbourhood of $X_0$, so one might try to construct the generic
fiber $X_\eta$ of $X_\infty$. As it turns out, this is indeed
possible, but we have to leave the category of (formal) schemes:
this generic fiber $X_\eta$ is a rigid variety over $K$.

\subsection{The affine case}
Let $A$ be an algebra topologically of finite type over $R$, and
consider the affine formal scheme $X_\infty=\mathrm{Spf}\,A$. The
tensor product $A\otimes_R K$ is a $K$-affinoid algebra, and the
generic fiber $X_\eta$ of $X_\infty$ is simply the $K$-affinoid
space
 $\Sp A\otimes_R K$.

Let $K'$ be any finite extension of $K$, and denote by $R'$ the
normalization of $R$ in $K'$. There exists a canonical bijection
between the set of morphisms of formal $R$-schemes
$\mathrm{Spf}\,R'\rightarrow X_\infty$, and the set of morphisms
of rigid $K$-varieties $\mathrm{Sp}\,K'\rightarrow X_\eta$.
Consider a morphism of formal $R$-schemes
$\mathrm{Spf}\,R'\rightarrow X_\infty$, or, equivalently, a
morphism of $R$-algebras $A\rightarrow R'$. Tensoring with $K$
yields a morphism of $K$-algebras $A\otimes_R K\rightarrow
R'\otimes_R K\cong K'$, and hence a $K'$-point of $X_\eta$.
Conversely, for any morphism of $K$-algebras $A\otimes_R
K\rightarrow K'$, the image of $A$ will be contained in $R'$,
since we've already seen in Section \ref{subsec-tate} that the
image of $R\{x_1,\ldots,x_m\}$ under any morphism of $K$-algebras
$T_m\rightarrow K^{alg}$ is contained in the normalization
$R^{alg}$ of $R$ in $K^{alg}$.

To any $R'$-section on $X_\infty$, we can associate a point of
$X_\infty$, namely the image of the singleton
$|\mathrm{Spf}\,R'|$. In this way, we obtain a specialization map
of sets
$$sp:|X_\eta|\rightarrow |X_\infty|=|X_0|$$
\subsection{The general case}\label{gengen}
The construction of the generic fiber for general $stft$ formal
$R$-schemes $X_\infty$ is obtained by gluing the constructions on
affine charts. The important point here is that the specialization
map $sp$ is continuous: if $X_\infty=\Spf A$ is affine, then for
any open formal subscheme $U_\infty$ of $X_\infty$, the inverse
image $sp^{-1}(U_\infty)$ is an admissible open in $X_\eta$; in
fact, if $U_\infty=\Spf B$ is affine, then $sp^{-1}(U_\infty)$ is
an affinoid domain in $X_\eta$, canonically isomorphic to
$U_\eta=\Sp B\otimes_R K$. Hence, the generic fibers of the
members $U_\infty^{(i)}$ of an affine open cover of a $stft$
formal $R$-scheme $X_\infty$ can be glued along the generic fibers
of the intersections $U_\infty^{(i)}\cap U_\infty^{(j)}$ to obtain
a rigid $K$-variety $X_\infty$, and the specialization maps glue
to a continuous map
$$sp:|X_\eta|\rightarrow |X_\infty|=|X_0|$$
This map can be enhanced to a morphism of ringed sites by
considering the unique morphism of sheaves
$$sp^{\sharp}:\mathcal{O}_{X_\infty}\rightarrow sp_*\mathcal{O}_{X_\eta}$$
which is given by the natural map
$$\mathcal{O}_{X_\infty}(U_\infty)=A\rightarrow A\otimes_R K=sp_*\mathcal{O}_{X_\eta}(U_\infty)$$ on any affine open
formal subscheme $U_\infty=\Spf A$ of $X_\infty$.

The generic fiber of a $stft$ formal $R$-scheme is a separated,
quasi-compact rigid $K$-variety. The formal scheme $X_\infty$ is
called a formal $R$-model for the rigid $K$-variety $X_\eta$.
Since the generic fiber is obtained by inverting $\pi$, it is
clear that the generic fiber does not change if we replace
$X_\infty$ by its maximal flat closed formal subscheme (by killing
$\pi$-torsion). If $K'$ is a finite extension of $K$ and $R'$ the
normalization of $R$ in $K'$, then we still have a canonical
bijection $X_\infty(R')=X_\eta(K')$.

The construction of the generic fiber is functorial: a morphism of
$stft$ formal $R$-schemes $h:Y_\infty\rightarrow X_\infty$ induces
a morphism of rigid $K$-varieties $h_{\eta}:Y_\eta\rightarrow
X_\eta$, and the square

$$\begin{CD}
Y_\eta@>h_\eta>> X_\eta \\@VspVV @VVspV \\Y_\infty @>h>> X_\infty
\end{CD}$$

commutes. We get a functor
$$(\,.\,)_\eta:(stft-For/R)\rightarrow (sqc-Rig/K):X_\infty\mapsto X_\eta$$
from the category of $stft$ formal $R$-schemes, to the category of
separated, quasi-compact rigid $K$-varieties.

For any locally closed subset $Z$ of $X_0$, the inverse image
$sp^{-1}(Z)$ is an admissible open in $X_\infty$, called the
\textit{tube} of $Z$ in $X_\infty$, and denoted by $]Z[$. If $Z$
is open in $X_0$, then $]Z[$ is canonically isomorphic to the
generic fiber of the open formal subscheme
$Z_\infty=(|Z|,\mathcal{O}_{X_\infty}|_{Z})$ of $X_\infty$. The
tube $]Z[$ is quasi-compact if $Z$ is open, but not in general.

Berthelot showed in \cite[0.2.6]{bert} how to construct the
generic fiber of a broader class of formal $R$-schemes, not
necessarily $tft$. If $Z$ is closed in $X_0$, then $]Z[$ is
canonically isomorphic to the generic fiber of the formal
completion of $X_\infty$ along $Z$ (this formal completion is the
locally topologically ringed space with underlying topological
space $|Z|$ and strucure sheaf
$$\lim_{\stackrel{\longleftarrow}{n}}\mathcal{O}_{X_\infty}/\mathcal{I}_Z^n$$
where $\mathcal{I}_Z$ is the defining ideal sheaf of $Z$ in
$X_\infty$). In particular, if $Z$ is a closed point $x$ of $X_0$,
then $]x[$ is the generic fiber of the formal spectrum of the
completed local ring $\widehat{\mathcal{O}}_{X_\infty,x}$ with its
adic topology (we did not define this notion; see
\cite[10.1]{ega1}).

\begin{example}\label{ex-tue}
Let $X_\infty$ be an affine $stft$ formal $R$-scheme, say
$X_\infty=\Spf A$. Consider tuple of elements $f_1,\ldots,f_r$ in
$A$, and denote by $Z$ the closed subscheme of $X_0$ defined by
the residue classes $\overline{f}_1,\ldots,\overline{f}_r$ in
$A_0$. The tube $]Z[$ of $Z$ in $X_\infty$  consists of the points
$x$ of $X_\eta$ with $|f_i(x)|<1$ for $i=1,\ldots,r$ (since this
condition is equivalent to $f_i(x)\equiv 0\mod (K^{alg})^{oo}$).

If $X_\infty=\Spf R\{x\}$, then $X_\eta$ is the closed unit disc
$\Sp K\{x\}$, and the special fiber $X_0$ is the affine line
$\A^1_k$. If we denote by $O$ the origin in $X_0$ and by $V$ its
complement, then $]V[$ is the affinoid domain $$U_1=\Sp
K\{x,T\}/(xT-1)$$ from Example \ref{strongtop} (the ``boundary''
of the closed unit disc\footnote{The notion of boundary is only
well-defined if you specify a center of the disc, since any point
of a closed disc can serve as a center, due to the ultrametric
property of the absolute value.}), and $]O[$ is the open unit disc
$U_2$ from the same example. The first one is quasi-compact, the
second is not.
\end{example}
\subsection{Localization by formal blow-ups}
The functor $(\,.\,)_\eta$ is not an equivalence. One can show
that formal blow-ups are turned into isomorphisms
\cite[4.1]{formrigI}. Intuitively, this is clear: the center
$\mathcal{I}$ of a formal blow-up contains a power of $\pi$, so it
becomes the unit ideal after inverting $\pi$.

In some sense, this is the only obstruction. Denote by
$\mathcal{C}$ the category of flat $stft$ formal $R$-schemes,
localized with respect to the formal blow-ups. This means that we
artificially add inverse morphisms for formal blow-ups, thus
turning them into isomorphisms.
The objects of $\mathcal{C}$ are simply the flat $stft$ formal
$R$-schemes, but a morphism in $\mathcal{C}$ from $Y_\infty$ to
$X_\infty$ is given by a triple $(Y'_\infty,\varphi_1,\varphi_2)$
where $\varphi_1:Y'_\infty\rightarrow Y_\infty$ is a formal
blow-up, and
 $\varphi_2:Y'_\infty\rightarrow
X_\infty$ a morphism of $stft$ formal $R$-schemes. We identify
this triple with another triple $(Y''_\infty,\psi_1,\psi_2)$ if
there exist a third triple $(Z_\infty,\chi_1,\chi_2)$ and
morphisms of $stft$ formal $R$-schemes $Z_\infty\rightarrow
Y'_\infty$ and $Z_\infty\rightarrow Y''_\infty$ such that the
obvious triangles commute.

Since admissible blow-ups are turned into isomorphisms by the
functor $(\,.\,)_\eta$, it factors though a functor
$\mathcal{C}\rightarrow (sqc-Rig/K)$. Raynaud \cite{Raynaud}
showed that this is an equivalence of categories (a detailed proof
is given in \cite{formrigI}). This means that the category of
separated, quasi-compact rigid $K$-varieties, can be described
entirely in terms of formal schemes. To give an idea of this
dictionary between formal schemes and rigid varieties, we list
some results. Let $X$ be a separated, quasi-compact rigid variety
over $K$.

\begin{itemize}
\item \cite[4.1(e),4.7]{formrigI} There exists a flat $stft$
formal $R$-scheme $X_\infty$ such that $X$ is isomorphic to
$X_\eta$.

 \item \cite[4.1(c+d)]{formrigI} If $X_\infty$ and
$Y_\infty$ are $stft$ formal $R$-schemes and
$\varphi:Y_\eta\rightarrow X_\eta$ is a morphism of rigid
$K$-varieties, then in general, $\varphi$ will not extend to a
morphism $Y_\infty\rightarrow X_\infty$ on the $R$-models.
However, by Raynaud's result, there exist a formal blow-up
$f:Y'_\infty\rightarrow Y_\infty$ and a morphism of $stft$ formal
$R$-schemes $g:Y'_\infty\rightarrow X_\infty$, such that
$\varphi=g_\eta\circ (f_\eta)^{-1}$. If $\varphi$ is an
isomorphism, we can find $(Y'_\infty,f,g)$ with both $f$ and $g$
formal blow-ups.

 \item \cite[4.4]{formrigI} For any
affinoid cover $\mathfrak{U}$ of $X$, there exist a formal model
$X_\infty$ of $X$ and a Zariski cover $\{U_1,\ldots,U_s\}$ of
$X_0$ such that $\mathfrak{U}=\{\,]U_1[,\ldots,]U_s[\,\}$.
\end{itemize}
See \cite{formrigI,formrigII,formrigIII,formrigIV} for many other
results.

\begin{example}
Consider the $stft$ formal $R$-schemes
\begin{eqnarray*}
X_\infty&=&\Spf R\{x\}/(x^2-1) \\Y_\infty&=&\Spf
R\{x\}/(x^2-\pi^2)\end{eqnarray*} The generic fibers $Y_\eta$ and
$X_\eta$ are isomorphic (both consist of two points $\Sp K$) but
it is clear that there is no morphism of $stft$ formal $R$-schemes
$Y_\infty\rightarrow X_\infty$ which induces an isomorphism
between the generic fibers. The problem is that the section
$x/\pi$ is not defined on $Y_\infty$; however, blowing up the
ideal $(x,\pi)$ adds this section to the ring of regular
functions, and the formal blow-up scheme is isomorphic to
$X_\infty$.

Next, consider the $stft$ formal $R$-scheme $Z_\infty=\Spf
R\{x\}$, and the standard cover of $Z_\eta$ defined by the couple
$(x,\pi)$. The cover consists of the closed disc $D(0,|\pi|)$ and
the closed annulus $Z_\eta\setminus D^{-}(0,|\pi|)$ (notation as
in Example \ref{ex-affinoid}). These sets are not tubes in
$Z_\infty$, since by Example \ref{ex-tue}, both sets have
non-empty intersection with the tube $]O[$ but do not coincide
with it. But if we take the formal blow-up $Z'_\infty\rightarrow
Z_\infty$ at the ideal $(x,\pi)$, then the rational subsets in our
standard cover are precisely the generic fibers of the blow-up
charts $\Spf R\{x,T\}/(xT-\pi)$ and $\Spf R\{x,T\}/(x-\pi T)$.
\end{example}

\subsection{Proper $R$-varieties}\label{subsec-proper}
Now let $X$ be a separated scheme of finite type over $R$, and
denote by $X_K$ its generic fiber. We denote by $(X_K)^o$ the set
of closed points of $X_K$. By \cite[0.3.5]{bert}, there exists a
canonical open immersion $\alpha:(\widehat{X})_{\eta}\rightarrow
(X_K)^{an}$. If $X$ is proper over $R$, then $\alpha$ is an
isomorphism.

For a proper $R$-scheme of finite type $R$, we can describe the
specialization map
$$sp:(X_K)^o=|(X_K)^{an}|=|\widehat{X}_\eta|\rightarrow |\widehat{X}|=|X_0|$$
as follows: let $x$ be a closed point of $X_K$, denote by $K'$ its
residue field, and by $R'$ the normalization of $R$ in $K'$. The
point $x$ defines a morphism $x:\Spec K'\rightarrow X$. The
valuative criterion for properness guarantees that the morphism
$\mathrm{Spec}\,R'\rightarrow \mathrm{Spec}\,R$ lifts to a unique
morphism $h:\mathrm{Spec}\,R'\rightarrow X$ with $h_K=x$. If we
denote by $0$ the closed point of $\mathrm{Spec}\,R'$, then
$sp(x)=h(0)\in |X_0|$.

In general, the open immersion
$\alpha:\widehat{X}_{\eta}\rightarrow (X_K)^{an}$ is strict.
Consider, for instance, a proper $R$-variety $X$, and let $X'$ be
the variety obtained by removing a closed point $x$ from the
special fiber $X_0$. Then $X'_K=X_K$; however, by taking the
formal completion $\widehat{X'}$, we loose all the points in
$\widehat{X}_\eta$ that map to $x$ under $sp$, i.e.
$\widehat{X'}_\eta=\widehat{X}_\eta\setminus \,]x[$. We'll see an
explicit example in the following section. This is another
instance of the fact that the rigid generic fiber
$\widehat{X'}_\eta$ is ``closer'' to the special fiber than the
scheme-wise generic fiber $X'_K$.

\subsection{Example: the projective line}
Let $X$ be the affine line $\Spec R[x]$ over $R$; then $X_K=\Spec
K[x]$, and $(X_K)^{an}$ is the rigid affine line $(\A^1_K)^{an}$
from Example \ref{ex-rigline}. On the other hand,
$\widehat{X}=\Spf R\{x\}$ and $\widehat{X}_\eta$ is the closed
unit disc $\Sp K\{x\}$. The canonical open immersion
$\widehat{X}_\eta\rightarrow (X_K)^{an}$ is an isomorphism onto
the affinoid domain in $(X_K)^{an}$ consisting of the points $z$
with $|x(z)|\leq 1$.

If we remove the origin $O$ from $X$, we get a scheme $X'$ with
$X'_K=X_K$. However, the formal completion of $X'$ is
$$\widehat{X'}=\Spf R\{x,T\}/(xT-1)$$ and its generic fiber is the
complement of $]O[$ in $\widehat{X}_\eta$ (see Example
\ref{ex-tue}).

Now let us turn to the projective line
$\mathbb{P}^1_R=\mathrm{Proj}\,R[x,y]$. The analytic projective
line $(\mathbb{P}^1_K)^{an}$ can be realized in different ways.
First, consider the usual affine cover of $\Proj^1_K$ by the
charts $U_1=\Spec K[x/y]$ and $U_2=\Spec K[y/x]$. Their
analytifications $(U_1)^{an}$ and $(U_2)^{an}$ are infinite unions
of closed discs (see Example \ref{ex-rigline}) centered at $0$,
resp. $\infty$. Gluing along the admissible opens
$(U_1)^{an}-\{0\}$ and $(U_2)^{an}-\{\infty\}$ in the obvious way,
we obtain $(\mathbb{P}^1_K)^{an}$.

On the other hand, we can look at the formal completion
$\widehat{\Proj^1_R}$. By the results in Section
\ref{subsec-proper}, we know that its generic fiber is canonically
isomorphic to $(\Proj^1_K)^{an}$. The $stft$ formal $R$-scheme
$\widehat{\Proj^1_R}$ is covered by the affine charts $V_1=\Spf
R\{x/y\}$ and $V_2=\Spf R\{y/x\}$ whose intersection is given by
$$V_0=\Spf R\{x/y,y/x\}/((x/y)(y/x)-1)$$ We've seen in Example
\ref{ex-tue} that the generic fibers of $V_1$ and $V_2$ are closed
unit discs around $x/y=0$, resp. $y/x=0$, and that $(V_0)_\eta$
coincides with their boundaries. So in this way,
$(\Proj^1_K)^{an}$ is realized as the Riemann sphere obtained by
gluing two closed unit discs along their boundaries.

\section{Berkovich spaces}\label{berkovich}
We recall some definitions from Berkovich' theory of analytic
spaces over non-archimedean fields. We refer to \cite{berk1}, or
to \cite{Berk2} for a short introduction. A very nice survey of
the theory and some of its applications is given in \cite{ducros}.

For a commutative Banach ring with unity $(\cA,\|\cdots \|)$, the
\textit{spectrum} $\cM(\cA)$ is the set of all bounded
multiplicative semi-norms $x:\cA\rightarrow \R_{+}$ (where
``bounded'' means that there exists a number $C>0$ such that
$x(a)\leq C\|a\|$ for all $a$ in $\cA$). If $x$ is a point of
$\cM(\cA)$, then $x^{-1}(0)$ is a prime ideal of $\cA$, and $x$
descends to an absolute value $|.|$ on the quotient field of
$\cA/x^{-1}(0)$. The completion of this field is called the
residue field of $x$, and denoted by $\cH(x)$. Hence, any point
$x$ of $\cM(\cA)$ gives rise to a bounded ring morphism $\chi_x$
from $\cA$ to the complete valued field $\cH(x)$, and $x$ is
completely determined by $\chi_x$. In this way, one can
characterize the points of $\cM(\cA)$ as equivalence classes of
bounded ring morphisms from $\cA$ to a complete valued field
\cite[1.2.2(ii)]{berk1}, just as one can view elements of the
spectrum $\Spec B$ of a commutative ring $B$ either as prime
ideals in $B$ or as equivalence classes of ring morphisms from $B$
to a field.

If we denote the image of $f\in \cA$ under $\chi_x$ by $f(x)$,
then $x(f)=|f(x)|$. We endow $\cM(\cA)$ with the weakest topology
such that $\cM(\cA)\rightarrow \R:x\mapsto |f(x)|$ is continuous
for each $f$ in $\cA$. This topology is called the
\textit{spectral topology} on $\cM(\cA)$. If $\cA$ is not the zero
ring, it makes $\cM(\cA)$ into a non-empty compact Hausdorff
topological space \cite[1.2.1]{berk1}. A bounded morphism of
Banach algebras $\cA\rightarrow \mathcal{B}$  induces a continuous
map $\cM(\mathcal{B})\rightarrow \cM(\cA)$ between their spectra.
In particular, the spectrum of $\cA$ only depends on the
equivalence class of $\|.\|$.

If $L$ is a non-archimedean field with non-trivial absolute value
and $A$ is an $L$-affinoid algebra (these are called
\textit{strictly} $L$-affinoid in Berkovich' theory) then $A$
carries a Banach norm, well-defined up to equivalence (Section
\ref{subsec-tate}). The spectrum $\cM(A)$ of $A$ is called a
strictly $L$-affinoid analytic space; Berkovich endows these
topological spaces with a structure sheaf of analytic functions.
General strictly $L$-analytic spaces are obtained by gluing
strictly $L$-affinoid spaces.

Any maximal ideal $x$ of $A$ defines a point of $\cM(A)$: the
bounded multiplicative semi-norm sending $f\in A$ to $|f(x)|$.
This defines a natural injection $\Sp A\rightarrow \cM(A)$, whose
image consists of the points $y$ of $\cM(A)$ with
$[\cH(y):L]<\infty$. So $\cM(A)$ contains the ``classical'' rigid
points of $\Sp A$, but in general also additional points $z$ with
$z^{-1}(0)$ not a maximal ideal. Beware that the natural map
$$\cM(A)\rightarrow \Spec A:z\mapsto z^{-1}(0)$$ is not injective,
in general: if $P\in \Spec A$ is not a maximal ideal, there may be
several bounded absolute values on $A/P$ extending the absolute
value on $L$.


For a Hausdorff strictly $L$-analytic space $X$, the set of rigid
points
$$X_{rig} :=\{x\in X\,|\,[\cH(x):L]<\infty\}$$ can be endowed
with the structure of a quasi-separated rigid variety over $L$ in
a natural way. Moreover, the functor $X\mapsto X_{rig}$ induces an
equivalence between the category of paracompact strictly
$L$-analytic spaces, and the category of quasi-separated rigid
varieties over $L$ which have an admissible affinoid covering of
finite type \cite[1.6.1]{Berk-etale}. The space $X_{rig}$ is
quasi-compact if and only if $X$ is compact.


The big advantage of Berkovich spaces is that they carry a
``true'' topology instead of a Grothendieck topology, with very
nice features (Hausdorff, locally connected by arcs,\ldots). As
we've seen, Berkovich obtains his spaces by adding points to the
 points of a rigid variety (not unlike the generic points in
algebraic geometry) which have an interpretation in terms of
valuations. We refer to \cite[1.4.4]{berk1} for a description of
the points and the topology of the closed unit disc
$D=\cM(\,L\{x\})$.

 To give a taste of these Berkovich spaces, let
us explain how two points of $D_{rig}$ can be joined by a path in
$D$. We assume, for simplicity, that $L$ is algebraically closed.
For each point $a$ of $D_{rig}$ and each $\rho\in [0,1]$ we define
$D(a,\rho)$ as the set of points $z$ in $D_{rig}$ with
$|x(z)-x(a)|\leq \rho$. This is not an affinoid domain if
$\rho\notin |L^*|$. Any such disc $E=D(a,\rho)$ defines a bounded
multiplicative semi-norm $|.|_{E}$ on the Banach algebra $L\{x\}$,
by mapping $f=\sum_{n=0}^{\infty}a_n(T-a)^n$ to
$$|f|_E=\sup_{z\in E}|f(z)|=\max_n |a_n|\rho^n$$ and hence, $E$ defines a
Berkovich-\textit{point} of $D$.
%
Now a path between two points $a,\,b$ of $D_{rig}=L^o$ can be
constructed as follows: put $\delta=|x(a)-x(b)|$ and consider the
path
$$\gamma:[0,1]\rightarrow D:t\mapsto \left\{\begin{array}{l}
                D(a,2t\delta),\ \textrm{if}\ 0\leq t\leq 1/2,
                \\ D(b,2(1-t)\delta),\ \textrm{if}\ 1/2\leq t\leq
                1.
                \end{array}\right.$$
Geometrically, this path can be seen as a closed disc around $a$,
 growing continuously in time $t$ until it contains $b$, and then
shrinking to $b$.

A remarkable feature of Berkovich' theory is that it can also be
applied to the case where $L$ carries the trivial absolute value.
If $k$ is any field, and $X$ is an algebraic variety over $k$,
then we can endow $k$ with the trivial absolute value and consider
the Berkovich analytic space $X^{an}$ associated to $X$ over $k$
\cite[3.5]{berk1}. Surprisingly, the topology of $X^{an}$ contains
non-trivial information on $X$. For instance, if $k=\C$, then the
rational singular cohomology $H_{sing}(X^{an},\Q)$ of $X^{an}$ is
canonically isomorphic to the weight zero part of the rational
singular cohomology of the complex analytic space $X(\C)$
\cite[1.1(c)]{berk-tate}. We refer to \cite{poineau} and
\cite{thuillier} for other applications of analytic spaces w.r.t.
trivial absolute values.

 Let us mention that there are still alternative
approaches to non-archimedean geometry, such as Fujiwara and
Kato's Zariski-Riemann spaces \cite{Fu-Ka}, or Huber's adic spaces
\cite{huber-adic}. See \cite{vdput-points} for a (partial)
comparison.
\section{Some applications}\label{applic}
\subsection{Relation to arc schemes and the Milnor fibration}\label{arcs}
\subsubsection{Arc spaces}
Let $k$ be any field, and let $X$ be a separated scheme of finite
type over $k$. Put $R=k[[t]]$. For each $n\geq 1$, we define a
functor
$$F_n:(k-alg)\rightarrow (Sets):A\rightarrow X(A\otimes_k R_n)$$
from the category of $k$-algebras to the category of sets. It is
representable by a separated $k$-scheme of finite type
$\mathcal{L}_n(X)$ (this is nothing but the Weil restriction of
$X\times_k R_n$ to $k$). For any pair of integers $m\geq n\geq 0$,
the truncation map $R_m\rightarrow R_n$ induces by Yoneda's Lemma
a morphism of $k$-schemes
$$\pi^m_n:\mathcal{L}_m(X)\rightarrow \mathcal{L}_n(X)$$
It is easily seen that these morphisms are affine, and hence, we
can consider the projective limit
$$\mathcal{L}(X):=\lim_{\stackrel{\longleftarrow}{n}}\mathcal{L}_n(X)$$
in the category of $k$-schemes. This scheme is called the arc
scheme of $X$. It satisfies $\mathcal{L}(X)(k')=X(k'[[t]])$ for
any field $k'$ over $k$ (these points are called arcs on $X$), and
comes with natural projections
$$\pi_n:\mathcal{L}(X)\rightarrow \mathcal{L}_n(X)$$
In particular, we have a morphism $\pi_0:\mathcal{L}(X)\rightarrow
\mathcal{L}_0(X)=X$. For any subscheme $Z$ of $X$, we put
$\mathcal{L}(X)_Z=\mathcal{L}(X)\times_X Z$. By Yoneda's lemma, a
morphism of separated $k$-schemes of finite type $h:Y\rightarrow
X$ induces $k$-morphisms $h:\mathcal{L}_n(Y)\rightarrow
\mathcal{L}_n(X)$, and by passage to the limit, a $k$-morphism
$h:\mathcal{L}(Y)\rightarrow \mathcal{L}(X)$.

If $X$ is smooth over $k$, the schemes $\mathcal{L}_n(X)$ and
$\mathcal{L}(X)$ are fairly well understood: if $X$ has pure
dimension $d$, then, for each pair of integers $m\geq n\geq 0$,
$\pi^m_n$ is a locally trivial fibration with fiber
$\A^{d(m-n)}_k$ (w.r.t. the Zariski topology). If $x$ is a
singular point of $X$, however, the scheme $\mathcal{L}(X)_x$ is
still quite mysterious. It contains a lot of information on the
singular germ $(X,x)$; interesting invariants can be extracted by
the theory of motivic integration (see \cite{DLinvent,DL3,Veys}).

The schemes $\mathcal{L}(X)_x$ and $\mathcal{L}(X)$ are not
Noetherian, in general, which complicates the study of their
geometric poperties. Already the fact that they have only finitely
many irreducible components if $k$ has characteristic zero, is a
non-trivial result. We will show how rigid geometry allows to
translate questions concerning the arc space into arithmetic
problems on rigid varieties.
\subsubsection{The
relative case}\label{relcase} Let $k$ be any algebraically closed
field of characteristic zero\footnote{This condition is only
imposed to simplify the arguments.}, and put $R=k[[t]]$. For each
integer $d>0$, $K=k((t))$ has a unique extension $K(d)$ of degree
$d$ in a fixed algebraic closure $K^{alg}$ of $K$, obtained by
joining a $d$-th root of $t$ to $K$. We denote by $R(d)$ the
normalization of $R$ in $K(d)$. For each $d>0$, we choose of a
$d$-th root of $t$ in $K^{alg}$, denoted by $\sqrt[d]{t}$, such
that $(\sqrt[de]{t})^e=\sqrt[d]{t}$ for each $d,e>0$. This choice
defines an isomorphism of $k$-algebras $R(d)\cong
k[[\sqrt[d]{t}]]$. It also induces an isomorphism of $R$-algebras
$$\varphi_d:R(d)\rightarrow R(d)':\sum_{i\geq 0}a_i(\sqrt[d]{t})^i\mapsto \sum_{i\geq 0}a_i t^i$$ where $R(d)'$
is the ring $R$ with $R$-algebra structure given by $$R\rightarrow
R:\sum_{i\geq 0}b_it^i\mapsto \sum_{i\geq 0}b_it^{id}$$

Let $X$ be a smooth irreducible variety over $k$, endowed with a
dominant morphism $f:X\rightarrow \A^1_k=\mathrm{Spec}\,k[t]$. We
denote by $\widehat{X}$ the formal completion of the $R$-scheme
$X_R=X\times_{k[t]} k[[t]]$; we will also call this the $t$-adic
completion of $f$. Its special fiber $X_0$ is simply the fiber of
$f$ over the origin.

There exists a tight connection between the points on the generic
fiber $\widehat{X}_\eta$ of $\widehat{X}$, and the arcs on $X$.
For any integer $d>0$, we denote by $\mathcal{X}(d)$ the closed
subscheme of $\mathcal{L}(X)$ defined by
$$\mathcal{X}(d)=\{\psi\in\mathcal{L}(X)\,|\,f(\psi)=t^d\}$$

We will construct a canonical bijection
$$\varphi:\widehat{X}_{\eta}(K(d))\rightarrow \mathcal{X}(d)(k)$$
such that the square

$$\begin{CD}
\widehat{X}_{\eta}(K(d))@>\varphi>> \mathcal{X}(d)(k)
\\ @VspVV @VV\pi_0V
\\ X_0(k)@>=>> X_0(k)
\end{CD}$$
commutes.

As we've seen in Section \ref{gengen}, the specialization morphism
of ringed sites $sp:\widehat{X}_\eta\rightarrow \widehat{X}$
induces a bijection $\widehat{X}_\eta(K(d))\rightarrow
\widehat{X}(R(d))$, and the morphism $sp$ maps a point of
$\widehat{X}_\eta(K(d))$ to the reduction modulo $\sqrt[d]{t}$ of
the corresponding point of $\widehat{X}(R(d))$.
 By Grothendieck's
Existence Theorem (Section \ref{completion}), the completion
functor induces a bijection $(X_R)(R(d))\rightarrow
\widehat{X}(R(d))$. Finally, the $R$-isomorphism
$\varphi_d:R(d)\rightarrow R(d)'$ induces a bijection
$$(X_R)(R(d))\rightarrow (X_R)(R(d)')=\mathcal{X}(d)(k)$$

In other words, if we take an arc
$\psi:\mathrm{Spec}\,R\rightarrow X$ with $f(\psi)=t^d$, then the
morphism $\widehat{\psi}_\eta$ yields a $K(d)$-point on $X_\eta$,
and this correspondence defines a bijection between
$\mathcal{X}(d)(k)$ and $\widehat{X}_\eta(K(d))$. Moreover, the
image of $\psi$ under the projection
$\pi_0:\mathcal{L}(X)\rightarrow X$ is nothing but the image of
the corresponding element of $\widehat{X}_\eta(K(d))$ under the
specialization map $sp:|\widehat{X}_\eta|\rightarrow
|\widehat{X}|=|X_0|$.

The Galois group $G(K(d)/K)=\mu_d(k)$ acts on
$\widehat{X}_\eta(K(d))$, and its action on the level of arcs is
easy to describe: if $\psi$ is an arc $\mathrm{Spec}\,R\rightarrow
X$ with $f(\psi)=t^d$, and $\xi$ is an element of $\mu_d(k)$, then
$\xi.\psi(t)=\psi(\xi.t)$.

The spaces $\mathcal{X}(d)$, with their $\mu_d(k)$-action, are
quite close to the arc spaces appearing in the definition of the
motivic zeta function associated to $f$ \cite[3.2]{DL3}. In fact,
the motivic zeta function can be realized in terms of the motivic
integral of a Gelfand-Leray form on $\widehat{X}_\eta$, and the
relationship between arc schemes and rigid varieties can be used
in the study of motivic zeta functions and the monodromy
conjecture, as is explained in \cite{NiSe3,NiSe}.
\subsubsection{The absolute case} This case is easily reduced to
the previous one. Let $X$ be any separated $k$-scheme of finite
type, and consider its base change $X_R=X\times_k R$. We denote by
$\widehat{X}$ the formal completion of $X_R$.

There exists a canonical bijection between the sets
$\mathcal{L}(X)(k)$ and $X_R(R)$. Hence, by the results in the
previous section, $k$-rational arcs on $X$ correspond to
$K$-points on the generic fiber $\widehat{X}_\eta$ of
$\widehat{X}$, by a canonical bijection
$$\varphi: \mathcal{L}(X)(k)\rightarrow \widehat{X}_\eta(K)$$
and the square

$$\begin{CD}
\mathcal{L}(X)(k)@>\varphi>> \widehat{X}_\eta(K)
\\ @V\pi_0VV @VVspV
\\ X(k) @>=>> X(k)
 \end{CD}$$
commutes. So the rigid counterpart of the space
$\mathcal{L}(X)_Z:=\mathcal{L}(X)\times_X Z$ of arcs with origin
in some closed subscheme $Z$ of $X$, is the tube $]Z[$ of $Z$ in
$\widehat{X}_\eta$ (or rather, its set of $K$-rational points).

Of course, the scheme structure on $\mathcal{L}(X)$ is very
different from the analytic structure on $\widehat{X}_\eta$.
Nevertheless, the structure on $\widehat{X}_\eta$ seems to be much
richer than the one on $\mathcal{L}(X)$, and one might hope that
some essential properties of the non-Noetherian scheme
$\mathcal{L}(X)$ are captured by the more ``geometric'' object
$\widehat{X}_\eta$. Moreover, there exists a satisfactory theory
of \'etale cohomology for rigid $K$-varieties (see for instance
\cite{Berk-etale} or \cite{huber-adic}), making it possible to
apply cohomological techniques to the study of the arc space.
\subsubsection{The analytic Milnor fiber}
Let $g:\C^m\rightarrow \C$ be an analytic map, and denote by $Y_0$
the analytic space defined by $g=0$. Let $x$ be a point of $Y_0$.
Consider an open disc $D:=B(0,\eta)$ of radius $\eta$ around the
origin in $\C$, and an open disc $B:=B(x,\varepsilon)$ in $\C^m$.
We denote by $D^{\times}$ the punctured disc $D-\{0\}$, and we put
$$X':=B\cap g^{-1}(D^{\times})$$ Then, for $0<\eta \ll
\varepsilon \ll 1$, the induced map
$$g':X'\longrightarrow
D^{\times}$$ is a $\mathcal{C}^{\infty}$ locally trivial
fibration, called the Milnor fibration of $g$ at $x$. It is
trivial if $g$ is smooth at $x$. Its fiber at a point $t$ of
$D^{\times}$ is denoted by $F_{x}(t)$, and it is called the
(topological) Milnor fiber of $g$ at $x$ (w.r.t. $t$). To remove
the dependency on the base point, one constructs the canonical
Milnor fiber $F_x$ by considering the fiber product
$$F_x:=X'\times_{D^{\times}}\widetilde{D^{\times}}$$ where
$\widetilde{D^{\times}}$ is the universal covering space
$$\widetilde{D^{\times}}=\{z\in \C\,|\,\Im(z)>-\log
\eta\}\rightarrow D^{\times}:z\mapsto \exp(iz)$$ Since this
covering space is contractible, $F_x$ is homotopically equivalent
to $F_x(t)$. The group of covering transformations
$\pi_1(D^{\times})$ acts on the singular cohomology of $F_x$; the
action of the canonical generator $z\mapsto z+2\pi$ of
$\pi_1(D^{\times})$ is called the monodromy transformation of $g$
at $x$. The Milnor fibration $g'$ was devised in \cite{Milnor} as
a tool to gather information on the topology of $Y_0$ near $x$.

We return to the algebraic setting: let $k$ be an algebraically
closed field of characteristic zero, put $R=k[[t]]$, let $X$ be a
smooth irreducible variety over $k$, and let $f:X\rightarrow
\A^1_k=\mathrm{Spec}\,k[t]$ be a dominant morphism. As before, we
denote by $\widehat{X}$ the formal $t$-adic completion of $f$,
with generic fiber $\widehat{X}_\eta$. For any closed point $x$ on
$X_0$, we put $\mathcal{F}_x:=]x[$, and we call this rigid
$K$-variety the analytic Milnor fiber of $f$ at $x$. This object
was introduced and studied in \cite{NiSe-Milnor,NiSe}. We consider
it as a bridge between the topological Milnor fibration and arc
spaces; a tight connection between these data is predicted by the
motivic monodromy conjecture. See \cite{NiSe3} for more on this
point of view.

The topological intuition behind the construction is the
following: the formal neighbourhood Spf\,$R$ of the origin in
$\A^1_k=\mathrm{Spec}\,k[t]$, corresponds to an infinitesimally
small disc around the origin in $\C$. Its inverse image under $f$
is realized as the $t$-adic completion of the morphism $f$; the
formal scheme $\widehat{X}$ should be seen as a tubular
neighborhood of the special fiber $X_0$ defined by $f$ on $X$. The
inverse image of the punctured disc becomes the ``complement'' of
$X_0$ in $\widehat{X}$, i.e. the generic fiber $\widehat{X}_\eta$
of $\widehat{X}$. The specialization map $sp$ can be seen as a
canonical ``contraction'' of $\widehat{X}_\eta$ on $X_0$, such
that $\mathcal{F}_x$ corresponds to the topological space $X'$
considered above. Note that this is not really the Milnor fiber
yet: we had to base-change to a universal cover of $D^{\times}$,
which corresponds to considering $\mathcal{F}_x\widehat{\times}_K
\widehat{K^{alg}}$ instead of $\mathcal{F}_x$, by the dictionary
between finite covers of $D^{\times}$ and finite extensions of
$K$. The monodromy action is translated into the Galois action of
$G(K^{alg}/K)\cong \widehat{\Z}(1)(k)$ on
$\mathcal{F}_x\widehat{\times}_K \widehat{K^{alg}}$.

It follows from the results in Section \ref{relcase} that, for any
integer $d>0$, the points in $\mathcal{F}_x(K(d))$ correspond
canonically to the arcs
$$\psi:\mathrm{Spec}\,k[[t]]\rightarrow X$$ satisfying
$f(\psi)=t^d$ and $\pi_0(\psi)=x$. Moreover, by Berkovich'
comparison result in \cite[3.5]{berk-vanish2} (see also Section
\ref{nearby}), there are canonical isomorphisms
$$H_{\acute{e}t}^i(\mathcal{F}_x\widehat{\times}_K
\widehat{K^{alg}},\Q_\ell)\cong R^i\psi_\eta(\Q_\ell)_x$$ such
that the Galois action of $G(K^{alg}/K)$ on the left hand side
corresponds to the monodromy action of $G(K^{alg}/K)$ on the
right. Here $H^*_{\acute{e}t}$ is \'etale $\ell$-adic cohomology,
and $R\psi_\eta$ denotes the $\ell$-adic nearby cycle functor
associated to $f$. In particular, if $k=\C$, this implies that
$H_{\acute{e}t}^i(\mathcal{F}_x\widehat{\times}_K
\widehat{K^{alg}},\Q_\ell)$ is canonically isomorphic to the
singular cohomology $H^i_{sing}(F_x,\Q_\ell)$ of the canonical
Milnor fiber $F_x$ of $f$ at $x$, and that the action of the
canonical topological generator of
$G(K^{alg}/K)=\widehat{\Z}(1)(\C)$ corresponds to the monodromy
trasformation, by Deligne's classical comparison theorem for
\'etale and analytic nearby cycles \cite[XIV]{sga7b}. In view of
the motivic monodromy conjecture, it is quite intriguing that
$\mathscr{F}_x$ relates certain arc spaces to monodromy action;
see \cite{NiSe3} for more background on this perspective.
\subsection{Deformation theory and lifting problems}
Suppose that $R$ has mixed characteristic, and let $X_0$ be a
scheme of finite type over the residue field $k$. Illusie sketches
in \cite[5.1]{formillusie} the following problem: is there a flat
scheme $X$ of finite type over $R$ such that $X_0=X\times_R k$?
Grothendieck suggested the following approach: first, try to
construct an inductive system $X_n$ of flat $R_n$-schemes of
finite type such that $X_n\cong X_m\times_{R_m}R_n$ for $m\geq
n\geq 0$. In many situations, the obstructions to lifting $X_n$ to
$X_{n+1}$ live in a certain cohomology group of $X_0$, and when
these obstructions vanish, the isomorphism classes of possible
$X_{n+1}$ correspond to elements in another approriate cohomology
group of $X_0$. Once we found such an inductive system, its direct
limit is a flat formal $R$-scheme $X_\infty$, topologically of
finite type. Next, we need to know if this formal scheme is
algebrizable, i.e. if there exists an $R$-scheme $X$ whose formal
completion $\widehat{X}$ is isomorphic to $X_\infty$. This scheme
$X$ would be a solution to our lifting problem. A useful criterion
to prove the existence of $X$ is the one quoted in Section
\ref{completion}: if $X_0$ is proper and carries an ample line
bundle that lifts to a line bundle on $X_\infty$, then $X_\infty$
is algebrizable. Moreover, the algebraic model $X$ is unique op to
isomorphism by Grothendieck's existence theorem (Section
\ref{completion}). For more concrete applications of this
approach, we refer to Section 5 of \cite{formillusie}.
\subsection{Nearby cycles for formal schemes}\label{nearby}
Berkovich used his \'etale cohomology theory for non-archimedean
analytic spaces, developed in \cite{Berk-etale}, to construct
nearby and vanishing cycles functors for formal schemes
\cite{Berk-vanish,berk-vanish2}. His formalism applies, in
particular, to $stft$ formal $R$-schemes $X_\infty$, and to formal
completions of such formal schemes along closed subschemes of the
special fiber $X_0$. Let us denote by $R\psi_\eta$ the functor of
nearby cycles, both in the algebraic and in the formal setting.
Suppose that $k$ is algebraically closed. Let $X$ be a variety
over $R$, and denote by $\widehat{X}$ its formal completion, with
generic fiber $\widehat{X}_\eta$. Let $Y$ be a closed subscheme of
$X_0$, and let $\mathcal{F}$ be an \'etale constructible sheaf of
abelian groups on $X\times_R K$, with torsion orders prime to the
characteristic exponent of $k$. Then Berkovich associates to
$\mathcal{F}$ in a canonical way an \'etale sheaf
$\widehat{\mathcal{F}}$ on $\widehat{X}_\eta$, and an \'etale
sheaf $\widehat{\mathcal{F}/Y}$ on the tube $]Y[$. His comparison
theorem \cite[3.1]{berk-vanish2} states that there are canonical
quasi-isomorphisms
$$R\psi_\eta(\mathcal{F})\cong
R\psi_\eta(\widehat{\mathcal{F}})\mbox{ and
}R\psi_\eta(\mathcal{F})|_{Y}\cong
R\psi_\eta(\widehat{\mathcal{F}/Y})$$ Moreover, by
\cite[3.5]{berk-vanish2} there is a canonical quasi-isomorphism
$$R\Gamma(Y,R\psi_\eta(\mathcal{F})|_{Y})\cong
R\Gamma(]Y[\widehat{\times}_K
\widehat{K^s},\widehat{\mathcal{F}/Y})$$ In particular, if $x$ is
a closed point of $X_0$, then $R^i\psi_\eta(\Q_\ell)_x$ is
canonically isomorphic to the $i$-th $\ell$-adic cohomology space
of the tube $]x[\widehat{\times}_K \widehat{K^s}$. Similar results
hold for tame nearby cycles and vanishing cycles.

This proves a conjecture of Deligne's, stating that
$R\psi_\eta(\mathcal{F})|_{Y}$ only depends on the formal
completion of $X$ along $Y$. In particular, the stalk of
$R\psi_\eta(\mathcal{F})$ at a closed point $x$ of $X_0$ only
depends on the completed local ring $\widehat{\mathcal{O}}_{X,x}$.
\subsection{Semi-stable reduction for curves}
Bosch and L\"utkebohmert show in
\cite{bosch-stable,bosch-stableII} how rigid geometry can be used
to construct stable models for smooth projective curves over a
non-archimedean field $L$, and uniformizations for Abelian
varieties. Let us briefly sketch their approach to stable
reduction of curves.

If $A$ is a reduced Tate algebra over $L$, then we define
\begin{eqnarray*}
A^{o}&=&\{f\in A\,|\,\|f\|_{sup}\leq 1\} \\A^{oo}&=&\{f\in
A\,|\,\|f\|_{sup}< 1\}
\end{eqnarray*}
Note that $A^o$ is a subring of $A$, and that $A^{oo}$ is an ideal
in $A^o$. The quotient $\widetilde{A}:=A^o/A^{oo}$ is a reduced
algebra of finite type over $\tilde{L}$, by
\cite[1.2.5.7+6.3.4.3]{BGR}, and
$\widetilde{X}:=\mathrm{Spec}\,\widetilde{A}$ is called the
canonical reduction of the affinoid space $X:=\mathrm{Sp}\,A$.
There is a natural reduction map $X\rightarrow \widetilde{X}$
mapping points of $X$ to closed points of $\widetilde{X}$. The
inverse image of a closed point $x$ of $\widetilde{X}$ is called
the formal fiber of $X$ at $x$; it is an open rigid subspace of
$X$.

Let $C$ be a projective, smooth, geometrically connected curve
over $L$, and consider its analytification $C^{an}$. By a
technical descent argument we may assume that $L$ is algebraically
closed. The idea is to construct a finite admissible cover
$\mathfrak{U}$ of $C^{an}$ by affinoid domains $U$ whose canonical
reductions $\widetilde{U}$ are semi-stable. If the cover $\mU$
satisfies a certain compatibility property, the canonical
reductions $\widetilde{U}$ can be glued to a semi-stable
$\widetilde{L}$-variety. From this cover $\mU$ one constructs a
semi-stable model for $C$. The advantage of passing to the rigid
world is that the Grothendieck topology on $C^{an}$ is much finer
then the Zariski topology on $C$, allowing finer patching
techniques.

To construct the cover $\mathfrak{U}$, it is proved that smooth
points and ordinary double points on $\widetilde{U}$ can be
recognized by looking at their formal fiber in $U$. For instance,
a closed point $x$ of $\widetilde{U}$ is smooth iff its formal
fiber is isomorphic to an open disc of radius $1$. An alternative
proof based on rigid geometry is given in \cite[5.6]{fresnel}.
\subsection{Constructing \'etale covers, and Abhyankar's Conjecture}
Formal and rigid patching techniques can also be used in the
construction of Galois covers; see \cite{patch} for an
introduction to this subject. This approach generalizes the
classical Riemann Existence Theorem for complex curves to a
broader class of base fields. Riemann's Existence Theorem states
that, for any smooth connected complex curve $X$, there is an
equivalence between the category of finite \'etale covers of $X$,
the category of finite analytic covering spaces of the complex
analytic space $X^{an}$, and the category of finite topological
covering spaces of $X(\C)$ (w.r.t. the complex topology). So the
problem of constructing an \'etale cover is reduced to the problem
of constructing a topological covering space, where we can proceed
locally w.r.t. the complex topology and glue the resulting local
covers. In particular, it can be shown in this way that any finite
group is the Galois group of a finite Galois extension of $\C(x)$,
by studying the ramified Galois covers of the complex projective
line.

The strategy in rigid geometry is quite similar: given a smooth
curve $X$ over a non-archimedean field $L$, we consider its
analytification $X^{an}$. We construct an \'etale cover $Y'$ of
$X^{an}$ by constructing covers locally, and gluing them to a
rigid variety. Then we use a GAGA-theorem to show that $Y'$ is
algebraic, i.e. $Y'=Y^{an}$ for some curve $Y$ over $L$; $Y$ is an
\'etale cover of $X$. Of course, several technical complications
have to be overcome to carry out this strategy.

We list some results that can be obtained by means of these
techniques, and references to their proofs.

\begin{itemize}
\item (Harbater) For any finite group $G$, there exists a ramified
Galois cover $f:X\rightarrow \mathbb{P}^1_L$ with Galois group
$G$, such that $X$ is absolutely irreducible, smooth, and
projective, and such that there exists a point $x$ in $X(L)$ at
which $f$ is unramified. An accessible proof by Q. Liu is given in
\cite{liu-harb}; see also \cite[\S\,3]{vdput-curves}. \item
(Abhyankar's Conjecture for the projective line) Let $k$ be an
algebraically closed field of characteristic $p>0$. A finite group
$G$ is the Galois group of a covering of $\mathbb{P}^1_k$, only
ramified over $\infty$, iff $G$ is generated by its elements
having order $p^n$ with $n\geq 1$. This conjecture was proven by
Raynaud in \cite{abhyankar-raynaud}. This article also contains an
introduction to rigid geometry and \'etale covers. \item
(Abhyankar's Conjecture) Let $k$ be an algebraically closed field
of characteristic $p>0$. Let $X$ be a smooth connected projective
curve over $k$ of genus $g$, let $\xi_0,\ldots,\xi_r$ ($r\geq 0$)
be distinct closed points on $X$, and let $\Gamma_{g,r}$ be the
topological fundamental group of a complex Riemann surface of
genus $g$ minus $r+1$ points (it is the free group on $2g+r$
generators). Put $U=X\setminus \{\xi_0,\ldots,\xi_r\}$. A finite
group $G$ is the Galois group of an unramified Galois cover of
$U$, iff every prime-to-$p$ quotient of $G$ is a quotient of
$\Gamma_{g,r}$. This conjecture was proven by Harbater in
\cite{abhyankar-harbater}.
\end{itemize}
\bibliographystyle{hplain}
\bibliography{wanbib,wanbib2}

\begin{thebibliography}{10}

\bibitem{sga7b}
{\em Groupes de monodromie en g\'eom\'etrie alg\'ebrique. {II}}.
\newblock Springer-Verlag, Berlin, 1973.
\newblock S\'eminaire de G\'eom\'etrie Alg\'ebrique du Bois-Marie 1967--1969
  (SGA 7 II), Dirig\'e par P. Deligne et N. Katz, Lecture Notes in Mathematics,
  Vol. 340.

\bibitem{berk1}
V.~G. Berkovich.
\newblock {\em {Spectral theory and analytic geometry over non-archimedean
  fields}}, volume~33 of {\em Mathematical Surveys and Monographs}.
\newblock AMS, 1990.

\bibitem{Berk-etale}
V.~G. Berkovich.
\newblock {\'Etale cohomology for non-Archimedean analytic spaces}.
\newblock {\em Publ. Math., Inst. Hautes \'Etud. Sci.}, 78:5--171, 1993.

\bibitem{Berk-vanish}
V.~G. Berkovich.
\newblock {Vanishing cycles for formal schemes}.
\newblock {\em Invent. Math.}, 115(3):539--571, 1994.

\bibitem{berk-vanish2}
V.~G. Berkovich.
\newblock {Vanishing cycles for formal schemes, {II}}.
\newblock {\em Invent. Math.}, 125(2):367--390, 1996.

\bibitem{Berk2}
V.~G. Berkovich.
\newblock {$p$-adic analytic spaces}.
\newblock {\em Doc. Math., J. DMV Extra Vol. ICM Berlin 1998, Vol. II}, pages
  141--151, 1998.

\bibitem{berk-tate}
V.~G. Berkovich.
\newblock {An analog of Tate's conjecture over local and finitely generated
  fields.}
\newblock {\em Int. Math. Res. Not.}, 2000(13):665--680, 2000.

\bibitem{bert}
P.~Berthelot.
\newblock {Cohomologie rigide et cohomologie rigide \`{a} supports propres}.
\newblock {\em Prepublication, Inst. Math. de Rennes}, 1996.

\bibitem{bosch}
S.~Bosch.
\newblock {\em {Lectures on formal and rigid geometry}}.
\newblock preprint, http://wwwmath1.uni-muenster.de/sfb/about/publ/bosch.html,
  2005.

\bibitem{BGR}
S.~Bosch, U.~G{\"{u}}ntzer, and R.~Remmert.
\newblock {\em {Non-Archimedean analysis. A systematic approach to rigid
  analytic geometry}}, volume 261 of {\em {Grundlehren der Mathematischen
  Wissenschaften}}.
\newblock Springer Verlag, 1984.

\bibitem{bosch-stableII}
S.~Bosch and W.~{L\"u}tkebohmert.
\newblock {Stable reduction and uniformization of abelian varieties. II}.
\newblock {\em Invent. Math.}, 78:257--297, 1984.

\bibitem{bosch-stable}
S.~Bosch and W.~{L\"u}tkebohmert.
\newblock {Stable reduction and uniformization of abelian varieties. I}.
\newblock {\em Math. Ann.}, 270:349--379, 1985.

\bibitem{formrigI}
S.~Bosch and W.~L{\"u}tkebohmert.
\newblock Formal and rigid geometry. {I}: {R}igid spaces.
\newblock {\em Math. Ann.}, 295(2):291--317, 1993.

\bibitem{formrigII}
S.~Bosch and W.~L{\"u}tkebohmert.
\newblock Formal and rigid geometry. {II}: {F}lattening techniques.
\newblock {\em Math. Ann.}, 296(3):403--429, 1993.

\bibitem{formrigIII}
S.~Bosch, W.~{L\"u}tkebohmert, and M.~Raynaud.
\newblock {Formal and rigid geometry. III: The relative maximum principle}.
\newblock {\em Math. Ann.}, 302(1):1--29, 1995.

\bibitem{formrigIV}
S.~Bosch, W.~{L\"u}tkebohmert, and M.~Raynaud.
\newblock {Formal and rigid geometry. IV: The reduced fibre theorem}.
\newblock {\em Invent. Math.}, 119(2):361--398, 1995.

\bibitem{DLinvent}
J.~Denef and F.~Loeser.
\newblock Germs of arcs on singular algebraic varieties and motivic
  integration.
\newblock {\em Invent. Math.}, 135:201--232, 1999, arxiv:math.AG/9803039.

\bibitem{DL3}
J.~Denef and F.~Loeser.
\newblock Geometry on arc spaces of algebraic varieties.
\newblock {\em Progr. Math.}, 201:327--348, 2001, arxiv:math.AG/0006050.

\bibitem{ducros}
A.~Ducros.
\newblock Espaces analytiques $p$-adiques au sens de {B}erkovich.
\newblock {\em Expos\'{e} 958 du s\'eminaire {B}ourbaki (mars 2006).}

\bibitem{fresnel}
J.~Fresnel and M.~van~der Put.
\newblock {\em {Rigid analytic geometry and its applications}}, volume 218 of
  {\em Progress in Mathematics}.
\newblock Boston, MA: Birk{h\"a}user, 2004.

\bibitem{Fu-Ka}
K.~Fujiwara and F.~Kato.
\newblock {Rigid geometry and applications}.
\newblock In {Mukai, S. et al.}, editor, {\em {Moduli spaces and arithmetic
  geometry}}, volume~45 of {\em {Advanced Studies in Pure Mathematics}}, pages
  327--386. {Mathematical Society of Japan, Tokyo}, 2006.

\bibitem{fga}
A.~Grothendieck.
\newblock {\em {Fondements de la g\'eom\'etrie alg\'ebrique. Extraits du
  S\'eminaire Bourbaki 1957-1962}}.
\newblock {Paris: Secr\'etariat math\'ematique }, 1962.

\bibitem{ega1}
A.~Grothendieck and J.~Dieudonn\'e.
\newblock El\'ements de {G}\'eom\'etrie {A}lg\'ebrique, {I}.
\newblock {\em Publ. Math., Inst. Hautes \'Etud. Sci.}, 4:5--228, 1960.

\bibitem{ega3}
A.~Grothendieck and J.~Dieudonn\'e.
\newblock El\'ements de {G}\'eom\'etrie {A}lg\'ebrique, {III}.
\newblock {\em Publ. Math., Inst. Hautes \'Etud. Sci.}, 11:5--167, 1961.

\bibitem{abhyankar-harbater}
D.~Harbater.
\newblock {Abhyankar's conjecture on Galois groups over curves}.
\newblock {\em Invent. Math.}, 117(1):1--25, 1994.

\bibitem{patch}
D.~Harbater.
\newblock {\em {Patching and Galois theory, \textrm{in} Schneps, L. (ed.),
  Galois groups and fundamental groups}}, volume~41 of {\em Math. Sci. Res.
  Inst. Publ.}, pages 313--424.
\newblock Cambridge: Cambridge University Press, 2003.

\bibitem{huber-adic}
R.~Huber.
\newblock {\em {\'Etale cohomology of rigid analytic varieties and adic
  spaces.}}, volume E30 of {\em {Aspects of Mathematics}}.
\newblock {Vieweg, Wiesbaden}, 1996.

\bibitem{formillusie}
L.~Illusie.
\newblock Grothendieck's {E}xistence {T}heorem in formal geometry.
\newblock In B.~Fantechi, L.~Goettsche, L.~Illusie, S.~Kleiman, N.~Nitsure, and
  A.~Vistoli, editors, {\em Fundamental Algebraic Geometry, Grothendieck's FGA
  explained}, Mathematical Surveys and Monographs. AMS, 2005.

\bibitem{liu-harb}
Q.~Liu.
\newblock Tout groupe fini est un groupe de galois sur $\mathbb{Q}_p(t)$,
  d'apr\`{e}s {H}arbater.
\newblock In M.~Fried, editor, {\em Recent Developments in the Inverse Galois
  Problem}, volume 186 of {\em Contemp. Math.}, pages 261--265. Providence, RI:
  AMS, 1995.

\bibitem{liu-rigid}
Q.~Liu.
\newblock {Une mini introduction \`a la g\'eom\'etrie analytique rigide}.
\newblock In B.~Deschamps, editor, {\em {Arithm\'etique des rev\^etements
  alg\'ebriques - Actes du colloque de Saint-\'Etienne }}, volume~5 of {\em
  S\'eminaires et Congr\`es}, pages 43--61. SMF, Paris, 2001.

\bibitem{formrig0}
W.~L{\"u}tkebohmert.
\newblock {Formal-algebraic and rigid-analytic geometry}.
\newblock {\em Math. Ann.}, 286(1-3):341--371, 1990.

\bibitem{Milnor}
J.~Milnor.
\newblock {\em Singular points of complex hypersurfaces}, volume~61 of {\em
  Annals of Math. Studies}.
\newblock Princeton University Press, 1968.

\bibitem{NiSe-Milnor}
J.~Nicaise and J.~Sebag.
\newblock Invariant de {S}erre et fibre de {M}ilnor analytique.
\newblock {\em C.R.Ac.Sci.}, 341(1):21--24, 2005.

\bibitem{NiSe}
J.~Nicaise and J.~Sebag.
\newblock The motivic {S}erre invariant, ramification, and the analytic
  {M}ilnor fiber.
\newblock {\em Invent. Math.}, 168(1):133--173, 2007.

\bibitem{NiSe3}
J.~Nicaise and J.~Sebag.
\newblock Rigid geometry and the monodromy conjecture.
\newblock In D.~Ch\'eniot et~al., editor, {\em Singularity Theory, Proceedings
  of the 2005 Marseille Singularity School and Conference}, pages 819--836.
  World Scientific, 2007.

\bibitem{poineau}
J.~Poineau.
\newblock {\em Espaces de Berkovich sur \Z}.
\newblock PhD thesis, Rennes, 2007.

\bibitem{Raynaud}
M.~Raynaud.
\newblock {G\'eom\'etrie analytique rigide d'apr\`es Tate, Kiehl, ... }.
\newblock {\em M\'emoires de la S.M.F.}, 39-40:319--327, 1974.

\bibitem{abhyankar-raynaud}
M.~Raynaud.
\newblock {Rev\^etements de la droite affine en caract\'eristique $p>0$ et
  conjecture d'Abhyankar}.
\newblock {\em Invent. Math.}, 116(1-3):425--462, 1994.

\bibitem{serre}
J.-P. Serre.
\newblock Classification des vari\'et\'es analytiques $p$-adiques compactes.
\newblock {\em Topology}, 3:409--412, 1965.

\bibitem{tate}
J.~Tate.
\newblock Rigid analytic geometry.
\newblock {\em Invent. Math.}, 12:257--289, 1971.

\bibitem{thuillier}
A.~Thuillier.
\newblock G\'eom\'etrie toro{\"{i}}dale et g\'eom\'etrie analytique non
  archim\'edienne. {A}pplication au type d'homotopie de certains sch\'emas
  formels.
\newblock {\em Manuscr. Math.}, 123(4):381--451, 2007.

\bibitem{vdput-curves}
M.~van~der Put.
\newblock {Valuation theory in rigid geometry and curves over valuation rings}.
\newblock In F.V. Kuhlmann, editor, {\em {Valuation theory and its
  applications. Volume I. Proceedings of the international conference and
  workshop, University of Saskatchewan, Saskatoon, Canada, July 28-August 11,
  1999}}, volume~32 of {\em Fields Inst. Commun.}, pages 369--383. Providence,
  RI: American Mathematical Society (AMS), 2002.

\bibitem{vdput-points}
M.~van~der Put and P.~Schneider.
\newblock {Points and topologies in rigid geometry}.
\newblock {\em Math. Ann.}, 302(1):81--103, 1995.

\bibitem{Veys}
W.~Veys.
\newblock Arc spaces, motivic integration and stringy invariants.
\newblock In S.~et~al. Izumiya, editor, {\em Singularity Theory and its
  applications}, volume~43 of {\em Advanced Studies in Pure Mathematics}, pages
  529--572. {Mathematical Society of Japan, Tokyo}, 2006.

\end{thebibliography}
\end{document}